\documentclass[12pt]{article}
\usepackage{latexsym,amssymb,amsmath}
\begin{document}

\baselineskip=19pt

\newcommand{\la}{\langle}
\newcommand{\ra}{\rangle}
\newcommand{\psp}{\vspace{0.4cm}}
\newcommand{\pse}{\vspace{0.2cm}}
\newcommand{\ptl}{\partial}
\newcommand{\dlt}{\delta}
\newcommand{\sgm}{\sigma}
\newcommand{\al}{\alpha}
\newcommand{\be}{\beta}
\newcommand{\G}{\Gamma}
\newcommand{\gm}{\gamma}
\newcommand{\vs}{\varsigma}
\newcommand{\lmd}{\lambda}
\newcommand{\td}{\tilde}
\newcommand{\vf}{\varphi}
\newcommand{\rd}{\mbox{Rad}}
\newcommand{\ad}{\mbox{ad}}
\newcommand{\stl}{\stackrel}
\newcommand{\ol}{\overline}
\newcommand{\ul}{\underline}
\newcommand{\es}{\epsilon}
\newcommand{\dmd}{\diamond}
\newcommand{\clt}{\clubsuit}
\newcommand{\vt}{\vartheta}
\newcommand{\ves}{\varepsilon}
\newcommand{\dg}{\dagger}
\newcommand{\kn}{\mbox{ker}}
\newcommand{\for}{\mbox{for}}
\newcommand{\dvg}{\mbox{div}}
\newcommand{\rar}{\rightarrow}
\newcommand{\NJ}{\mathbb{N}^{\ell_1+\ell_2}}
\newcommand{\ZJ}{\mathbb{Z}^{\ell_1+\ell_2}}
\newcommand{\bs}{\backslash}
\newcommand{\der}{\mbox{Der}\:}
\newcommand{\lra}{\Longleftrightarrow}
\newcommand{\BB}{\mathbb}

\begin{center}{\Large \bf Structure of Divergence-Free Lie Algebras}\end{center}
\vspace{0.2cm}

\begin{center}{\large Yucai Su$^{\ast}$ and  Xiaoping Xu$^{\dag}$}\end{center}

* Department of Applied Mathematics, Shanghai Jiaotong University, 1954 Huashan Road, Shanghai 200030, P. R. China.

\dag Department of Mathematics, The Hong Kong University of Science \& Technology, Clear Water Bay, Kowloon, Hong Kong, P. R. China.

\vspace{0.3cm}

\begin{abstract}{One of  the four well-known series of simple Lie algebras
of Cartan type is the series of Lie algebras of Special type, which are divergence-free Lie algebras associated with polynomial algebras and the operators of taking partial derivatives, connected with volume-preserving
diffeomorphisms. In this paper, we determine the structure space of the divergence-free Lie algebras associated with  pairs of a commutative associative algebra with an identity element and its finite-dimensional commutative locally-finite derivation subalgebra such that the commutative associative algebra is derivation-simple with respect to the derivation subalgebra.}\end{abstract}

\section{Introduction} 

Volume-preserving diffeomorphisms are important transformations on manifolds. Lie groups of volume-preserving diffeomorphisms are fundamental structures in the theory of idealized fluids. They have been studied through their Lie algebras, whose elements are divergence-free. One of  the four well-known series of simple Lie algebras of Cartan type is the series of Lie algebras of Special type, which are divergence-free Lie algebras associated with polynomial algebras and the operators of taking partial derivatives. Graded generalizations of the Lie algebras of Special type have been studied by Kac [K1], Osborn [O],
Djokovic and Zhao [DZ], and Zhao [Z]. Supersymmetric graded generalizations of the algebras have been investigated by Kac [K2], [K4]. Motivated from the nongradedness of the Lie algebras generated by conformal algebras (cf. [K3], [X2]) and the Lie algebras associated with vertex operator algebras (cf. [B], [FLM]), the second author [X1]  introduced a natural class of nongraded generalizations of the Lie algebras of Special type, which are simple algebras in general. The term ``divergence-free'' is used because we want to give the readers an intuitive picture of these Lie algebras.  In this paper, we shall determine the isomorphism classes of the divergence-free Lie algebras  introduced in [X2].
Below we shall give a detailed technical introduction. 

Throughout this paper, let $\BB{F}$ be a field  with characteristic 0. All the vector spaces are assumed over $\BB{F}$. Denote by $\BB{Z}$ the ring of integers and by $\BB{N}$ the set of nonnegative numbers $\{0,1,2,3,...\}$. We shall always identify $\BB{Z}$ with $\BB{Z}1_{\BB{F}}$ when the context is clear. 

Let ${\cal A}=\BB{F}[t_1,t_2,...,t_n]$ be the algebra of polynomials in $n$ variables.
 A {\it derivation} $\ptl$  of ${\cal A}$ is linear transformation of ${\cal A}$ such that
$$\ptl(uv)=\ptl(u)v+u\ptl(v)\qquad\for\;\;u,v\in{\cal A}.\eqno(1.1)$$
Typical derivations are $\{\ptl_{t_1},\ptl_{t_2},...,\ptl_{t_n}\}$,  the operators of taking partial derivatives. The space $\der{\cal A}$ of all the derivations of ${\cal A}$ forms a Lie algebra. Identifying the elements of 
${\cal A}$ with their corresponding multiplication operators, we have
$$\der{\cal A}=\sum_{i=1}^n{\cal A}\ptl_{t_i},\eqno(1.2)$$
which forms a simple Lie algebra. The Lie algebra $\der{\cal A}$ is called a {\it Witt algebra of rank} $n$, usually denoted as ${\cal W}(n,\BB{F})$. The Lie algebra ${\cal W}(n,\BB{F})$ acts on the Grassmann algebra $\hat{\cal A}$ of differential forms over ${\cal A}$ as follows.
$$\ptl(df)=d(\ptl(f)),\;\;\ptl(\omega\wedge \nu)=\ptl(\omega)\wedge\nu +\omega\wedge\ptl(\nu)\eqno(1.3)$$
for $f\in{\cal A},\;\omega,\nu\in\hat{\cal A},\;\ptl\in {\cal W}(n,\BB{F})$.

For $\ptl=\sum_{i=1}^nf_i\ptl_{t_i}\in {\cal W}(n,\BB{F})$, we define the {\it divergence of} $\ptl$ by
$$\dvg\:\ptl=\sum_{i=1}^n\ptl_{t_i}(f_i).\eqno(1.4)$$
Set 
$${\cal S}(n,\BB{F})=\{\ptl\in {\cal W}(n,\BB{F})\mid \dvg\:\ptl =0\}.\eqno(1.5)$$
Then ${\cal S}(n,\BB{F})$ is the divergence-free Lie subalgebra of ${\cal W}(n,\BB{F})$. Moreover, ${\cal S}(n,\BB{F})$ is the subalgebra of ${\cal W}(n,\BB{F})$ annihilating the volume form; that is,
$${\cal S}(n,\BB{F})=\{\ptl\in {\cal W}(n,\BB{F})\mid \ptl(dt_1\wedge dt_2\wedge\cdots\wedge dt_n)=0\}.\eqno(1.6)$$
When $\BB{F}=\BB{R}$ the field of real numbers, $e^{\ptl}$ is a volume-preserving diffeomorphism for each $\ptl\in {\cal S}(n,\BB{F})$. The algebra ${\cal S}(n,\BB{F})$ is called a {\it Lie algebra of Special type}.

For any positive integer $n$, an additive subgroup $G$ of $\BB{F}^n$ is called {\it nondegenerate} if $G$ contains an $\BB{F}$-basis of $\BB{F}^n$. Let $\ell_1,\;\ell_2$ and $\ell_3$ be three nonnegative integers such that
$$\ell=\ell_1+\ell_2+\ell_3>0.\eqno(1.7)$$
Take any nondegenerate additive subgroup $\G$ of $\BB{F}^{\ell_2+\ell_3}$ and $\G=\{0\}$ when $\ell_2+\ell_3=0$.  Let ${\cal A}(\ell_1,\ell_2,\ell_3;\G)$ be a free $\BB{F}[t_1,t_2,...,t_{\ell_1+\ell_2}]$-module with the basis
$$\{x^{\al}\mid \al\in\G\}.\eqno(1.8)$$
Viewing ${\cal A}(\ell_1,\ell_2,\ell_3;\G)$ as a vector space over $\BB{F}$, we define a commutative associative algebraic operation ``$\cdot$'' on ${\cal A}(\ell_1,\ell_2,\ell_3;\G)$ by
$$(\zeta x^{\al})\cdot (\eta x^{\be})=\zeta\eta x^{\al+\be}\qquad\for\;\;\zeta,\eta\in \BB{F}[t_1,t_2,...,t_{\ell_1+\ell_2}],\;\al,\be\in\G.\eqno(1.9)$$
Note that  $x^{0}$ is the identity element, which is denoted as $1$ for convenience. When the context is clear, we shall omit the notation ``$\cdot$'' in any associative algebra product.

We define the linear transformations 
$$\{\ptl_{t_1},...,\ptl_{t_{\ell_1+\ell_2}},\ptl^{\ast}_1,...,\ptl^{\ast}_{\ell_2+\ell_3}\}\eqno(1.10)$$
on ${\cal A}(\ell_1,\ell_2,\ell_3;\G)$ by 
$$\ptl_{t_i}(\zeta x^{\al})=\ptl_{t_i}(\zeta)x^{\al},\;\;\ptl^{\ast}_j(\zeta x^{\al})=\al_j \zeta x^{\al}\eqno(1.11)$$
for $\zeta\in \BB{F}[t_1,t_2,...,t_{\ell_1+\ell_2}]$ and $\al=(\al_1,...,\al_{\ell_2+\ell_3})\in\G.$
 Then $\{\ptl_{t_1},...,\ptl_{t_{\ell_1+\ell_2}},\ptl^{\ast}_1,...,\ptl^{\ast}_{\ell_2+\ell_3}\}$ are mutually commutative derivations of ${\cal A}(\ell_1,\ell_2,\ell_3;\G)$.

Throughout this paper, we shall use the following notation of index set
$$\ol{m,n}=\{m,m+1,m+2,...,n\}\qquad\for\;\;m,n\in\Bbb{Z},\;m\leq n.\eqno(1.12)$$
 Moreover, we also take $\ol{m,n}=\emptyset$ if $m>n$. Set
$$\ptl_i=\ptl_{t_i},\;\;\ptl_{\ell_1+j}=\ptl^{\ast}_j+\ptl_{t_{\ell_1+j}},\;\;\ptl_{\ell_1+\ell_2+l}=\ptl_{\ell_2+l}^{\ast}\eqno(1.13)$$
for $i\in\ol{1,\ell_1},\;j\in\ol{1,\ell_2}$ and $l\in\ol{1,\ell_3}$. Then $\{\ptl_i\mid i\in\ol{1,\ell}\}$ is an $\BB{F}$-linearly independent set of derivations. 
Let
$${\cal D}=\sum_{i=1}^{\ell}\BB{F}\ptl_i\eqno(1.14)$$
and 
$${\cal W}(\ell_1,\ell_2,\ell_3;\G)= {\cal A}(\ell_1,\ell_2,\ell_3;\G){\cal D}.\eqno(1.15)$$
Then ${\cal W}(\ell_1,\ell_2,\ell_3;\G)$ is a simple Lie algebra of Witt type constructed by the second author in [X1].
 
A linear transformation $T$ on a vector space $V$ is called {\it locally-finite} if
$$\dim(\mbox{Span}\:\{T^n(u)\mid n\in\Bbb{N}\})<\infty\qquad\for\;\;u\in V.\eqno(1.16)$$
 Zhang and the authors of this paper proved that the pairs $({\cal A}(\ell_1,\ell_2,\ell_3;\G),{\cal D})$ for different parameters $(\ell_1,\ell_2,\ell_3;\G)$ enumerate all the  pairs $({\cal A},{\cal D})$ of a commutative associative algebra ${\cal A}$ with an identity element and its finite-dimensional commutative locally-finite derivation subalgebra ${\cal D}$ such that the commutative associative algebra is derivation-simple with respect to ${\cal D}$ and
$$\bigcap_{d\in{\cal D}}\kn_{d}=\BB{F}.\eqno(1.17)$$

Denote by $M_{m\times n}$ the algebra of $m\times n$ matrices with entries in $\BB{F}$ and by $GL_m$ the group of invertible $m\times m$ matrices with entries in $\BB{F}$. Set
$$G_{\ell_2,\ell_3}=\left\{\left(\begin{array}{cc}A&0_{\ell_2\times\ell_3}\\B& C\end{array}\right)\mid A\in GL_{\ell_2},\;B\in M_{\ell_2\times\ell_3},\;C\in  GL_{\ell_3}\right\},\eqno(1.18)$$
where $0_{\ell_2\times\ell_3}$ is the $\ell_2\times\ell_3$ matrix whose entries are zero. Then $G_{\ell_2,\ell_3}$ forms a subgroup of $GL_{\ell_2+\ell_3}$. Define an action of $G_{\ell_2,\ell_3}$ on $\BB{F}^{\ell_2+\ell_3}$ by
$$ g(\al)=\al g^{-1}\;\;(\mbox{matrix multiplication})\qquad\for\;\;\al\in\BB{F}^{\ell_2+\ell_3},\;g\in G_{\ell_2,\ell_3}.\eqno(1.19)$$
For any nondegenerate additive subgroup $\Upsilon$ of $\BB{F}^{\ell_2+\ell_3}$ and $g\in G_{\ell_2,\ell_2}$, the set
$$g(\Upsilon)=\{g(\al)\mid \al\in \Upsilon\}\eqno(1.20)$$
also forms a nondegenerate additive subgroup of $\BB{F}^{\ell_2+\ell_3}$. Let
$$\Omega_{\ell_2+\ell_3}=\mbox{the set of nondegenerate additive subgroups of}\;\BB{F}^{\ell_2+\ell_3}.\eqno(1.21)$$
We have an action of $G_{\ell_2,\ell_3}$ on $\Omega_{\ell_2+\ell_3}$ by (1.20). Define the moduli space
$${\cal M}^W_{\ell_2,\ell_3}=\Omega_{\ell_2+\ell_3}/G_{\ell_2,\ell_3},\eqno(1.22)$$
which is the set of $G_{\ell_2,\ell_3}$-orbits in $\Omega_{\ell_2+\ell_3}$. 

Another related result in [SXZ] is as follows. The Lie algebras ${\cal W}(\ell_1,\ell_2,\ell_3;\G)$ and ${\cal W}(\ell_1',\ell_2',\ell_3';\G')$ are isomorphic if and only if $(\ell_1,\ell_2,\ell_3)=(\ell_1',\ell_2',\ell_3')$ and there exists an element $g\in  G_{\ell_2,\ell_3}$ such that
$g(\G)=\G'$. In particular, there exists a one-to-one correspondence between the set of isomorphism classes of the Lie algebras of the form (1.15) and the following set:
$$SW=\{(\ell_1,\ell_2,\ell_3,\varpi)\mid (0,0,0)\neq(\ell_1,\ell_2,\ell_3)\in\BB{N}^3,\;\varpi\in {\cal M}^W_{\ell_2,\ell_3}\}.\eqno(1.23)$$
In other words, the set $SW$  is the structure space of the simple Lie algebras of Witt type in the form (1.15). 

We define the divergence by
$$\dvg\:\ptl=\sum_{i=1}^{\ell}\ptl_i(u_i)\qquad\for\;\;\ptl=\sum_{i=1}^{\ell}u_i\ptl_i\in {\cal W}(\ell_1,\ell_2,\ell_3;\G),\eqno(1.24)$$
and set
$${\cal S}(\ell_1,\ell_2,\ell_3;\G)=\{\ptl \in {\cal W}(\ell_1,\ell_2,\ell_3;\G)\mid \dvg\:\ptl=0\}.\eqno(1.25)$$
Let $\rho\in\G$ be any element. Then the space
$${\cal S}(\ell_1,\ell_2,\ell_3;\rho,\G)=x^{\rho}{\cal S}(\ell_1,\ell_2,\ell_3;\G)\eqno(1.26)$$
forms a Lie subalgebra of the Lie algebra ${\cal W}(\ell_1,\ell_2,\ell_3;\G)$, which is simple if $\ell_1+\ell_2>0$ or $\rho=0$ or $\ell\geq 3$ by the proofs of the simplicity of the Lie algebras of Special type in [DZ] and [X1]. When $\ell_1+\ell_2=0,\;\ell_3=2$ and $\rho\neq 0$,  the derived subalgebra $({\cal S}(\ell_1,\ell_2,\ell_3;\rho,\G))^{(1)}$ is simple and has codimension one in ${\cal S}(\ell_1,\ell_2,\ell_3;\rho,\G)$.  The Lie algebra ${\cal S}(\ell_1,\ell_2,\ell_3;\rho,\G)$ was introduced  by the second author in [X1], as a nongraded generalization of graded simple Lie algebras of Special type. We also call it a {\it divergence-free Lie algebra} in this paper. The special case ${\cal S}(\ell_1,0,\ell_3;\rho,\G)$ was studied by Osborn [O],
Djokovic and Zhao [DZ], and Zhao [Z].

Define the moduli space
$${\cal M}^{\cal S}_{\ell_2,\ell_3}=(\G\times \Omega_{\ell_2+\ell_3})/G_{\ell_2,\ell_3},\eqno(1.27)$$
where the action of $G_{\ell_2,\ell_3}$ on $\G\times\Omega_{\ell_2+\ell_3}$
is defined by $g(\rho,\Upsilon)=(g(\rho),g(\Upsilon))$
for $(\rho,\Upsilon)\in\G\times\Omega_{\ell_2+\ell_3}$
(cf. (1.18)-(1.21)).
The main theorem of this paper is as follows.
\psp

{\bf Main Theorem}. {\it The Lie algebras} ${\cal S}(\ell_1,\ell_2,\ell_3;\rho,\G)$
{\it and} ${\cal S}(\ell_1',\ell_2',\ell_3';\rho',\G')$ {\it with} $\ell\geq 3$ {\it are isomorphic if and only if} $(\ell_1,\ell_2,\ell_3)=(\ell_1',\ell_2',\ell_3')$
{\it and there exists an element} $g\in  G_{\ell_2,\ell_3}$ {\it such that}
$g(\G)=\G'$, {\it and} $g(\rho)=\rho'$ {\it if} $\ell_1=0$. {\it In particular, there exists a one-to-one correspondence between the set of isomorphism classes of the Lie algebras of the form (1.26) and the set} $SW$ {\it in (1.23) if} $\ell_1>0$, {\it and between the set of isomorphism classes of the Lie algebras of the form (1.26) and  the following set}:
$$SS=\{(\ell_1,\ell_2,\ell_3,\varpi)\mid (0,0,0)\neq(\ell_1,\ell_2,\ell_3)\in\BB{N}^3,\;\varpi\in {\cal M}^{\cal S}_{\ell_2,\ell_3}\}\eqno(1.28)$$
{\it if} $\ell_1=0$.
\psp

When $\ell=2$, the Lie algebra ${\cal S}(\ell_1,\ell_2,\ell_3;\rho,\G)$  is also a Lie algebra of Hamiltonian type.  We shall determine the structure of the Lie algebras of Hamiltonian type in a subsequent paper. 

The paper is organized as follows. In Section 2, we shall present some basic properties of the divergence-free Lie algebras. Structure of derivation algebras of the divergence-free Lie algebras will be determined in Section 3. Using the results in the above two sections, we shall give the proof of the main theorem in Section 4.

\section{Some Basic Properties of the Lie algebras}

In this section, we shall present some basic properties of the divergence-free Lie algebra ${\cal S}(\ell_1,\ell_2,\ell_3;\rho,\G)$. The related assumptions and settings are the same as in the paragraphs of (1.7)-(1.15) and (1.24)-(1.26).

For convenience, we denote
$${\cal D}_1=\sum_{i=1}^{\ell_1}\BB{F}\ptl_i,\;\;{\cal D}_2=\sum_{j=1}^{\ell_2}\BB{F}\ptl_{\ell_1+j},\;\;{\cal D}_3=\sum_{p=1}^{\ell_3}\BB{F}\ptl_{\ell_1+\ell_2+p}\eqno(2.1)$$
(cf. (1.13)) and
$${\cal A}={\cal A}(\ell_1,\ell_2,\ell_3;\G),\;\;{\cal W}={\cal W}(\ell_1,\ell_2,\ell_3;\G),\;\;{\cal S}={\cal S}(\ell_1,\ell_2,\ell_3;\rho,\G)\eqno(2.2)$$
(cf. (1.9), (1.15), (1.26)). Then
$${\cal D}={\cal D}_1\oplus{\cal D}_2\oplus{\cal D}_3\eqno(2.3)$$
(cf.  (1.14)). Moreover, we define
$$D_{p,q}(u)=x^{\rho}(\ptl_q(x^{-\rho}u)\ptl_p-
\ptl_p(x^{-\rho}u)\ptl_q)\qquad \for\;\;p,q\in\ol{1,\ell},\;u\in {\cal A}.
\eqno(2.4)$$
It can be proved that
$${\cal S}=\mbox{Span}\:\{D_{p,q}(u)\mid p,q\in\ol{1,\ell};\;u\in {\cal A}\}.\eqno(2.5)$$

Denote
$$t^{\vec{i}}=t_1^{i_1}t_2^{i_2}\cdots t_{\ell_1+\ell_2}^{i_{\ell_1+\ell_2}},\;\;x^{\al,\vec{i}}=t^{\vec{i}}x^{\al}\qquad\for\;\;\vec{i}=(i_1,...,i_{\ell_1+\ell_2})\in
\NJ,\;\al\in\G\eqno(2.6)$$
and
$$a_{[i]}=(0,...,0,\stl{i}{a},0,...,0)\qquad\for\;\;a\in\BB{F}.\eqno(2.7)$$
In the rest of this paper, we shall use the convention that if a notation has not been defined but technically appears in an expression, we treat it as zero. For instance, when we use the notations $\al\in \G$ and $\vec i\in\NJ$, we teat $\al_p=0$ and $i_{\ell_1+\ell_2+q}=0$ for $p\leq 0$ and $q>0$ if they appear in an expression. Note
\begin{eqnarray*}D_{p,q}(x^{\al,\vec{i}})&=&x^{\al,\vec{i}}(((\al_{q-\ell_1}-\rho_{q-\ell_1})\ptl_p-(\al_{p-\ell_1}-\rho_{p-\ell_1})\ptl_q)\\ & &+i_qx^{\al,\vec{i}-1_{[q]}}\ptl_p-i_px^{\al,\vec{i}-1_{[p]}}\ptl_q\hspace{7.5cm}(2.8)\end{eqnarray*}
for $p,q\in\ol{1,\ell},\;\al\in\G$ and $\vec{i}\in\NJ$. Set
$$
{\cal S}_{\al}=\mbox{Span}\:\{D_{p,q}(x^{\al,\vec{i}})\mid p,q\in\ol{1,\ell};\;\vec{i}\in\NJ\}\qquad \for\;\;\al\in\G.\eqno(2.9)$$
Then
$${\cal S}=\bigoplus_{\al\in\G}{\cal S}_{\al},\eqno(2.10)$$
is a $\G$-graded Lie algebra, and
$${\cal S}_0=\mbox{Span}\:\{t^{\vec{i}}(-\rho_{q-\ell_1}\ptl_p+\rho_{p-\ell_1}\ptl_q)+i_qt^{\vec{i}-1_{[q]}}\ptl_p-i_pt^{\vec{i}-1_{[p]}}\ptl_q\mid p,q\in\ol{1,\ell},\;\vec{i}\in\NJ\}\eqno(2.11)$$
is a Lie subalgebra of ${\cal S}$.

We define a bilinear product from ${\cal D}\times \G\rar\BB{F}$ by
$$\al(\ptl)=\la \ptl,\al\ra=\sum_{i=1}^{\ell_2+\ell_3}a_{\ell_1+i}
\al_i\qquad\for\;\;\ptl=\sum_{i=1}^{\ell}a_i\ptl_i\in{\cal D},\;
\al=(\al_1,...,\al_{\ell_2+\ell_3})\in\G.\eqno(2.12)$$
For any $\al\in\G$, let
$${\cal D}_{\al}=\left\{\begin{array}{ll}
\{0\}&\mbox{ if}\;\ell_1=\ell_2=0,\al=0,\\ {\cal D}_3&\mbox{ if}\;\ell_1+\ell_2=1,\al=0,\\ \{\ptl\in{\cal D}\mid\la\ptl,\al\ra=0\}&\mbox{ otherwise.}\end{array}\right.\eqno(2.13)$$
In particular, we have
$$\dim {\cal D}_{\al}\geq\ell-1\geq 2\;\;\mbox{if}\;\;\ell_1+\ell_2\geq 1\;\;\mbox{or}\;\; \al\neq 0.\eqno(2.14)$$
Denote
$$ |\vec{i}|=\sum_{p=1}^{\ell_1+\ell_2}i_p\qquad\for\;\;\vec{i}\in\NJ.\eqno(2.15)$$
Define a total order on $\NJ$ by:
$$\vec{i}<\vec{j}\;\;\mbox{if}\;\;|\vec{i}|<|\vec{j}|\;\;\mbox{or}\;\;|\vec{i}|=|\vec{j}|\;\;
\mbox{with}\;\;i_1=j_1,...,i_{p-1}=j_{p-1},i_p<j_p\eqno(2.16)$$
for some $p\in\ol{1,\ell_1+\ell_2}$. For any
$$u=x^{\al}(t^{\vec{i}}\ptl_{\vec{i}}+\sum_{\vec{j}<\vec{i}}t^{\vec{j}}\ptl_{\vec{j}})\in {\cal S}_{\al}\;\;\mbox{with}\;\;\ptl_{\vec{i}},\ptl_{\vec{j}}\in{\cal D}\;\;\mbox{such that}\;\;\ptl_{\vec{i}}
\neq 0,\eqno(2.17)$$
 we define
$$u_{ld}=x^{\al,\vec{i}}\ptl_{\vec{i}}
\eqno(2.18)$$
to be the {\it leading term} of $u$. In this case, we say that $u$ has
{\it leading level} $|\vec{i}|$ and that $u$ has {\it leading degree} $\vec{i}$, denoted by $d(u_{ld})=\vec{i}$.

Set
$${\cal S}_{\al}^{[\vec{i}]}=\{u\in {\cal S}_{\al}\mid d(u_{ld})\le\vec{i}\},\ \
{\cal S}_{\al}^{(\vec{i})}=\{u\in {\cal S}_{\al}\mid d(u_{ld})<\vec{i}\}.\eqno(2.19)$$
and
$${\cal S}^{[\vec{i}]}=\bigoplus_{\al\in\G}{\cal S}_{\al}^{[\vec{i}]},\ \
{\cal S}^{(\vec{i})}=\bigoplus_{\al\in\G}{\cal S}_{\al}^{(\vec{i})}.\eqno(2.20)$$
Then 
$${\cal S}^{[0]}=\bigoplus_{\al\in\G}{\cal S}_{\al}^{[0]}\eqno(2.21)$$
 is a Lie subalgebra of ${\cal S}$.
\psp

{\bf Lemma 2.1}. {\it For} $\al\in\G$, {\it we have}
$${\cal S}_{\al}^{[0]}=\{x^{\al}\ptl\mid\ptl\in{\cal D}_{\al-\rho}\}.\eqno(2.22)$$

{\it Proof}. First we assume $\al\in\G\bs\{\rho\}$.
There exists $q\in\ol{\ell_1+1,\ell}$ such that $\al_{q-\ell_1}\neq\rho_{q-\ell_1}$.
For $\ptl=\sum_{p=1}^\ell a_p\ptl_p\in{\cal D}_{\al-\rho}$, we have
\begin{eqnarray*}\hspace{2cm}\sum_{p=1}^\ell a_p D_{p,q}(x^{\al})&=&x^{\al}(\sum_{p=1}^\ell a_p(\al_{q-\ell_1}-\rho_{q-\ell_1})\ptl_p+\sum_{p=1}^\ell a_p(\al_{p-\ell_1}-\rho_{p-\ell_1})\ptl_q)\\&=&x^{\al}((\al_{q-\ell_1}-\rho_{q-\ell_1})\ptl+\la\ptl,\al-\rho\ra\ptl_q)\\&=&(\al_{q-\ell_1}-\rho_{q-\ell_1})x^{\al}\ptl\hspace{6cm}(2.23)\end{eqnarray*}
by (2.13). Thus $x^{\al}\ptl\in {\cal S}_{\al}^{[0]}$.
Conversely, (2.8) and (2.9) imply
\begin{eqnarray*}\hspace{2cm}{\cal S}_{\al}^{[0]}&=&\mbox{Span}\:\{D_{p,q}(x^{\al,\vec{i}})\mid\vec{i}=\vec{0}\mbox{ or }\al_{p-\ell_1}-\rho_{p-\ell_1}=0,\\& &\al_{q-\ell_1}-\rho_{q-\ell_1}=0\mbox{ and }\vec{i}=1_{[p]} \mbox{ or }1_{[q]}
\},\hspace{4.4cm}(2.24)\end{eqnarray*}
and we see that ${\cal S}_{\al}^{[0]}\subseteq\{x^{\al}\ptl\mid \ptl\in{\cal D}_{\al-\rho}\}$.

Next we assume $\al=\rho$. If $\ell_1+\ell_2=0$, then  we have  ${\cal S}_{\rho}^{[0]}=\{0\}$ by (2.8) and (2.24). Suppose  $\ell_1+\ell_2\geq 1$. Let $p\in\ol{1,\ell}$. If $\ell_1+\ell_2\geq 2$ or
$p\in\ol{\ell_1+\ell_2+1,\ell}$,
then we can choose $q\in\ol{1,\ell_1+\ell_2}\bs\{p\}$ and get
$$x^{\rho}\ptl_p=D_{p,q}(x^{\rho,1_{[q]}})\in {\cal S}_{\rho}^{[0]}.\eqno(2.25)$$
However, if $\ell_1+\ell_2=1$ and $p=1$, then  (2.8) and (2.24) imply $x^{\rho}\ptl_1\notin {\cal S}_{\rho}^{[0]}$. So we obtain (2.22) by (2.13).$\qquad\Box$
\psp

We define a linear function $\chi_p:{\cal D}\rar\BB{F}$ by
$$\chi_p(\ptl)=a_p\qquad\for\;\;\ptl=\sum_{p=1}^\ell a_p\ptl\in{\cal D}.\eqno(2.26)$$
\psp

{\bf Lemma 2.2}. {\it For} $\al\in\G,\;\vec{i}\in\NJ$ {\it and} $\ptl
=\sum_{p=1}^\ell a_p\ptl_p\in{\cal D}\bs\{0\}$, {\it there exists an element in} ${\cal S}_{\al}$ {\it with leading term} $x^{\al,\vec{i}}\ptl$ {\it if and only if} $\ptl\in{\cal D}_{\al-\rho}$, {\it and}  $i_{\ell_1+\ell_2}=0$ {\it or} $a_{\ell_1+\ell_2}=0$ {\it if} $\al=\rho$.
\psp

{\it Proof}. Suppose $\al_{q-\ell_1}\neq\rho_{q-\ell_1}$ for some  $q\in\ol{\ell_1+1,\ell}$. Let $\ptl\in{\cal D}_{\al-\rho}$. As in
(2.23),  the element
\begin{eqnarray*}& &(\al_{q-\ell_1}-\rho_{q-\ell_1})^{-1}\sum_{p=1}^\ell a_p \ptl_{p,q}(x^{\al,\vec{i}})\\ &=&x^{\al,\vec{i}}\ptl+(\al_{q-\ell_1}-\rho_{q-\ell_1})^{-1}\sum_{p=1}^\ell a_p(i_q x^{\al,\vec{i}-1_{[q]}}\ptl_p-i_p x^{\al,\vec{i}-1_{[p]}}\ptl_q)\hspace{3.7cm}(2.27)\end{eqnarray*}
is in ${\cal S}_{\al}$ with the leading term $x^{\al,\vec{i}}\ptl$.
Conversely, if ${\cal S}_{\al}$ has an element with the leading term $x^{\al,\vec{i}}\ptl$,
 then (2.8) and (2.9) imply $\ptl\in{\cal D}_{\al-\rho}$.

 Assume $\al=\rho$. Let $\ptl=\sum_{p=1}^{\ell}a_p\ptl_p\in{\cal D}_0$.
If $a_{\ell_1+\ell_2}=0$, then 
\begin{eqnarray*}\hspace{2cm}& &(i_{\ell_1+\ell_2}+1)^{-1}\sum_{p=1}^\ell a_pD_{p,\ell_1+\ell_2}(x^{\rho,\vec{i}+1_{[\ell_1+\ell_2]}})\\&=&x^{\rho,\vec{i}}\ptl-(i_{\ell_1+\ell_2}+1)^{-1}\sum_{p=1}^{\ell_1+\ell_2-1}a_p i_p x^{\rho,\vec{i}-1_{[p]}+1_{[\ell_1+\ell_2]}}\ptl_{\ell_1+\ell_2}\hspace{2.7cm}(2.28)\end{eqnarray*}
is an element in ${\cal S}_{\rho}$ with the leading term $x^{\rho,\vec{i}}\ptl$. If $a_{\ell_1+\ell_2}\ne0$, then $i_{\ell_1+\ell_2}=0$ by our assumption. By (2.13), we must have $\ell_1+\ell_2\geq 2$. So we can write $\ptl=a_{\ell_1+\ell_2}\ptl_{\ell_1+\ell_2}+\ptl'$ such that $\ptl'\in{\cal D}_0$ and $\chi_{\ell_1+\ell_2}(\ptl')=0$. Then (2.28) shows that there
exists an element in ${\cal S}_{\rho}$ with the leading term $x^{\rho,\vec{i}}\ptl'$. Moreover, we have
$$x^{\rho,\vec{i}}\ptl_{\ell_1+\ell_2}=(i_1+1)^{-1}D_{\ell_1+\ell_2,1}(x^{\rho,\vec{i}+1_{[1]}})\in {\cal S}_{\rho}.\eqno(2.29)$$
Thus there exists an element in ${\cal S}_{\rho}$ with leading term $x^{\rho,\vec{i}}\ptl$. Similarly, if an element of ${\cal S}_{\rho}$ has a leading term $x^{\rho,\vec{i}}\ptl$, then we have $\ptl\in{\cal D}_0$ and $\chi_{\ell_1+\ell_2}(\ptl)=0$ or $i_{\ell_1+\ell_2}=0$ by (2.8) and (2.9).$\qquad\Box$
\psp

{\bf Lemma 2.3}. {\it If} $\ell_1\ge1$, then ${\cal S}(\ell_1,\ell_2,\ell_3;\rho,\G)\cong
{\cal S}(\ell_1,\ell_2,\ell_3;\rho',\G)$ {\it for any} $\rho,\rho'\in\G$.
\psp

{\it Proof}. For any $\al^{(1)},...,\al^{(\ell_1)}\in\G$, we define an automorphism $\psi$ of the associate algebra ${\cal A}$ (cf. (2.2)) by
$$
\psi(x^{\al,\vec i})=x^{\al+\sum_{p=1}^{\ell_1}i_p\al^{(p)},\vec i}\qquad\for\;\;(\al,\vec{i})\in\G\times\NJ.\eqno(2.30)$$
This induces a Lie algebra automorphism $\ol\psi$ on ${\cal W}$ (cf. (2.2)):
$$\ol{\psi}(x^{\al,\vec i}\ptl_q)=x^{\al+\sum_{p=1}^{\ell_1}i_p\al^{(p)}-\al^{(q)},\vec{i}}\ptl_q\qquad \for\;\;(\al,\vec{i})\in\G\times\NJ,\;q\in\ol{1,\ell_1};\eqno(2.31)$$
$$\ol{\psi}(x^{\al,\vec i}\ptl_q)=x^{\al+\sum_{p=1}^{\ell_1}i_p\al^{(p)},\vec{i}}(\ptl_q-\sum_{p=1}^{\ell_1}t_p\ptl_p)\qquad \for\;\;(\al,\vec{i})\in\G\times\NJ,\;q\in\ol{\ell_1+1,\ell}.\eqno(2.32)$$
It can be verified by (1.24)-(1.26) that the automorphism $\ol\psi$ of ${\cal W}$ induces an isomorphism
$$\ol\psi: {\cal S}(\ell_1,\ell_2,\ell_3;\rho,\G)\rar {\cal S}(\ell_1,\ell_2,\ell_3;\rho^{\ast},\G)\;\;\mbox{with}\;\;\rho^{\ast}=\rho+\sum_{p=1}^{\ell_1}\al^{(p)}.\eqno(2.33)$$
Taking $\al^{(1)}=\rho'-\rho$ and $\al^{(2)}=\cdots=\al^{(\ell_1)}=0$, we have $\rho^{\ast}=\rho'.\qquad\Box$
\psp

For convenience, we shall always take
$$\rho=0\qquad\mbox{if}\;\;\ell_1\ge1\eqno(2.34)$$
in the rest of this paper.

\section{Structure of the Derivation Algebras}

In this section, we shall determine the structure of the derivation algebra of the divergence-free Lie algebra ${\cal S}(\ell_1,\ell_2,\ell_3;\rho,\G)$ defined by (1.26), which has been simply denoted as ${\cal S}$.

In the rest of this paper, 
$$\mbox{we choose an element}\;\;\bar{x}^{\al,\vec{i}}\ptl\in{\cal S}\;\;\mbox{with the leading term}\;\;x^{\al,\vec{i}}\ptl
\eqno(3.1)$$
for each
$$(\al,\vec{i})\in \G\times\NJ,\;\ptl\in{\cal D}_{\al-\rho}\;\mbox{such that}\;i_{\ell_1+\ell_2}=0\;\mbox{or}\;\chi_{\ell_1+\ell_2}(\ptl)=0\;\mbox{if}\;\al=\rho.\eqno(3.2)$$

Recall that a derivation $d$ of the Lie algebra ${\cal S}$ is a linear transformation on ${\cal S}$ such that
$$d([u_1,u_2])=[d(u_1),u_2]+[u_1,d(u_2)]\qquad\for\;\;u_1,u_2\in {\cal S}.
\eqno(3.3)$$
 Denote 
$$\der{\cal S}=\mbox{the space of the derivations of}\;{\cal S}.\eqno(3.4)$$
It is well known that $\der{\cal S}$ forms a Lie algebra with respect to the commutator of linear transformations of ${\cal S}$, and $\ad_{\cal S}$ is an ideal of $\der{\cal S}$.

Recall that ${\cal S}$ is a $\G$-graded Lie algebra with the grading defined by (2.8) and (2.9). For $\al\in\G$, we set
$$(\der{\cal S})_{\al}=\{d\in\der{\cal S}\mid d({\cal S}_{\be})\subset{\cal S}_{\al+\be}\;\for\;\be\in\G\}.\eqno(3.5)$$

Fix an element $\al\in\G$. Recall the notations in (2.20). Consider a nonzero element
$$ d\in(\der{\cal S})_{\al}\;\;\mbox{such that}\;\;d({\cal S}^{[\vec j]})\subset{\cal S}^{[\vec i+\vec j]}\;\;\mbox{for any}\;\;\vec{j}\in\NJ,\eqno(3.6)$$
where $\vec{i}\in\ZJ$ is a fixed element.
Define the order in $\ZJ$ as in (2.16). We remark that
if $\vec j\in\ZJ\bs\NJ$, it does not mean that ${\cal S}^{[\vec j]}=0$ since
it is possible that the level $|\vec j|>0$. But if $\vec j\in\ZJ\bs\NJ$, then we have ${\cal S}^{[\vec j]}={\cal S}^{(\vec j)}$ (cf. (2.20)). Moreover, we can assume
$$\vec{i}\;\;\mbox{in (3.6) is the minimal element satisfying the condition}.\eqno(3.7)$$

For any $(\be,\vec{j})\in\G\times\NJ$, we define a linear map $e_{\be,\vec j}: {\cal D}_{\be-\rho}\rar{\cal D}_{\al+\be-\rho}$ by
$$d(\bar{x}^{\be,\vec j}\ptl)\equiv \bar{x}^{\al+\be,\vec i+\vec j}
e_{\be,\vec j}(\ptl)\;\; ({\rm mod\,}{\cal S}^{(\vec i+\vec j)}).\eqno(3.8)$$

Using the assumption (3.2) and the notation in (3.1), we have
$$[\bar{x}^{\be,\vec j}\ptl',\bar{x}^{\gm,\vec k}\ptl'']\equiv\bar{x}^{\be+\gm,\vec j+\vec k}
(\gm(\ptl')\ptl''-\be(\ptl'')\ptl')\ ({\rm mod\,}{\cal S}^{(\vec j+\vec k)}).
\eqno(3.9)$$
Applying $d$ to (3.9), we get
\begin{eqnarray*}& &\la e_{\be,\vec j}\ptl',\gm\ra\ptl''-\la\ptl'',\al+\be\ra e_{\be,\vec j}\ptl'
+\la\ptl',\al+\gm\ra e_{\gm,\vec k}\ptl''
-\la e_{\gm,\vec k}\ptl'',\be\ra\ptl'\\ &=&e_{\be+\gm,\vec j+\vec k}(\la\ptl',\gm\ra\ptl''-\la\ptl'',\be\ra\ptl').\hspace{8.6cm}(3.10)\end{eqnarray*}
\pse

{\bf Lemma 3.1}. {\it If} ${\cal D}_{\rho}\not\subset{\cal D}_{\al}$ {\it and} $\al\neq 0$, {\it then the derivation} $d$ {\it in (3.6) is an inner derivation of} ${\cal S}$.
\psp

{\it Proof}. By assumption, there exists $\hat{\ptl}\in{\cal D}_{\rho}\bs{\cal D}_{\al}$
with $\la\hat{\ptl},\al\ra=1$. Fix such a $\hat{\ptl}$ and set
$$\td{\ptl}=-e_{0,\vec{0}}\hat{\ptl}.\eqno(3.11)$$
Taking $\gm=0,\;\vec k=\vec{0}$ and $\ptl''=\hat{\ptl}\in{\cal D}_{\gm-\rho}$ in (3.10), we obtain
$$-\la\hat{\ptl},\al+\be\ra e_{\be,\vec j}\ptl'+\la\ptl',\al\ra e_{0,\vec{0}}\hat{\ptl}-\la e_{0,\vec{0}}\hat{\ptl},\be\ra\ptl'=-\la\hat{\ptl},\be\ra e_{\be,\vec j}\ptl'.\eqno(3.12)$$
By (3.11), (3.12) becomes
$$e_{\be,\vec j}\ptl'=\la\td{\ptl},\be\ra\ptl'-\la\ptl',\al\ra\td{\ptl}\qquad \mbox{for any}\;\;(\be,\vec{j},\ptl')\;\;\mbox{satisfying (3.2)}.\eqno(3.13)$$
Thus $e_{\be,\vec j}$ does not depend on $\vec j$ for any $\be$. Letting
$\vec j=\vec 0$ in (3.8), we have $\vec i\in\NJ$ by (3.7). Moreover, we have $\al\ne\rho$ by the assumption ${\cal D}_{\rho}\not\subset{\cal D}_{\al}$. Letting $\be=2\rho$ and
$\ptl'=\hat{\ptl}\in{\cal D}_{\be-\rho}$ in (3.13) and taking the bilinear product (2.12) of (3.13) with
$\be+\al-\rho$, we obtain $\la\td{\ptl},\al-\rho\ra=0$ because $\la e_{\be,\vec j}\hat{\ptl},\be+\al-\rho\ra=0$ by (3.8).  Set
$$d'=d-\ad_{\bar{x}^{\al,\vec{i}}\td{\ptl}}.\eqno(3.14)$$
By (3.8) and (3.13), we obtain that
$$d'(\bar{x}^{\be,\vec j}\ptl')\in {\cal S}^{(\vec i+\vec j)}\qquad \mbox{for any}\;\;(\be,\vec{j},\ptl')\;\;\mbox{satisfying (3.2)}.\eqno(3.15)$$
By induction on $\vec{i}$ in (3.7), $d'$ is an inner derivation. Hence $d$ is an inner derivation.$\qquad\Box$
\psp

Set
$${\cal W}^{[0]}_{\rho}=x^{\rho}{\cal D}\eqno(3.16)$$
(cf. (1.8), (1.14)).
\psp

{\bf Lemma 3.2}. {\it Let} $d,\;\vec i$ {\it be the same as in (3.6) and (3.7) with} $\al\neq 0$. {\it If} ${\cal D}_{\rho}\subseteq{\cal D}_{\al}$, {\it then} $\vec i\in\NJ$ {\it and}
$d$ {\it is an inner derivation when} $\al\ne\rho$, {\it or} $d\in\ad_{{\cal W}^{[0]}_{\rho}}|_{\cal S}+\ad_{\cal S}$ {\it  when} $\al=\rho$.
\psp

{\it Proof}. Notice that any element of $\ad_{{\cal W}^{[0]}_{\rho}}|_{\cal S}$ are homogeneous derivations of ${\cal S}$ of degree $\rho$. By (2.13), the condition ${\cal D}_{\rho}\subseteq{\cal D}_{\al}$ and the assumption (2.34) imply  $\ell_1=0$. So we have $\rho=a\al$ for some $a\in\BB{F}$. Note that in this case
$\ell_2+\ell_3=\ell\ge3$.
\pse

{\it Claim 1}. For any $\be\in\G\bs\BB{F}\al$ and
$\vec j\in\mathbb{N}^{\ell_1+\ell_2}$,
there exists
$a_{\be,\vec j}\in\BB{F}$ such that
$$e_{\be,\vec j}\ptl'=a_{\be,\vec j}\ptl'\qquad\for\;\; \ptl'\in{\cal D}_{\al}\cap{\cal D}_{\be}.
\eqno(3.17)$$

Let $\be$ and $\ptl'$ be the same as in the above. Taking $\gm=-\be+\rho$ and $\ptl''\in{\cal D}_{\gm-\rho}={\cal D}_{\be}$ in (3.10), we have that the third term of the left-hand side vanishes due to
$\ptl'\in{\cal D}_{\al}\cap{\cal D}_{\be}$, and the right-hand side also vanishes because of $\la\ptl',\gm\ra=\la\ptl'',\be\ra=0$. Thus (3.10) becomes
$$\la e_{\be,\vec j}\ptl',-\be+\rho\ra\ptl''-\la\ptl'',\al\ra e_{\be,\vec j}\ptl'-\la e_{-\be+\rho,\vec k}\ptl'',\be\ra\ptl'=0\eqno(3.18)$$
for $\ptl'\in{\cal D}_{\al}\cap{\cal D}_{\be}$ and $\ptl''\in{\cal D}_{\be}$. Taking the product (2.12) of (3.18) with $-\be+\rho$, we obtain
$$\la e_{\be,\vec j}\ptl',-\be+\rho\ra
\la\ptl'',\rho-\al\ra=0\eqno(3.19)$$
by the facts that  $\ptl'\in{\cal D}_{\al}\cap{\cal D}_{\be}$ and $\ptl''\in{\cal D}_{\be}$. If $\al\ne\rho$, then we take $\ptl''\in{\cal D}_{\be}\bs{\cal D}_{\al}$, and (3.19) forces $\la e_{\be,\vec j}\ptl',-\be+\rho\ra=0$. So (3.18) implies (3.17) with
$$a_{\be,\vec j}=-\la\ptl'',\al\ra^{-1}
\la e_{-\be+\rho,0},\ptl''\ra.\eqno(3.20)$$

Suppose  $\al=\rho$. Then $\rho\ne0$ by our assumption $\al\ne0$. We want to prove that the coefficient of the first term of
(3.18) is zero. Taking $\gm=-\rho$ and $\ptl''=\ptl'\in{\cal D}_{\rho}\cap{\cal D}_{\be}\bs\{0\}$ in (3.10), we have that
the second term and the third term of the left-hand side
vanish, and so does the right-hand side. Thus (3.10) is equivalent to
$$\la e_{\be,\vec j}\ptl', -\rho\ra=\la e_{-\rho,\vec k}\ptl',\be\ra\qquad\for\;\;\ptl'\in{\cal D}_{\rho}\cap{\cal D}_{\be}.\eqno(3.21)$$
Next we take $\gm=-\rho$, $\ptl''\in{\cal D}_{\rho}\cap{\cal D}_{\be}\bs\{0\}$ and $\ptl'\in{\cal D}_{\be-\rho}\bs
{\cal D}_{\rho}$. Then (3.10) becomes
$$\la e_{\be,\vec j}\ptl',-\rho\ra\ptl''-\la e_{-\rho,\vec k}\ptl'',\be\ra\ptl'=\la\ptl',-\rho\ra e_{\be-\rho,\vec j+\vec k}\ptl''.
\eqno(3.22)$$
Taking the product (2.12) of (3.22) with $\be$, we get
$$-\la e_{-\rho,\vec k}\ptl'',\be\ra\la\ptl',\be\ra=\la\ptl',-\rho\ra \la e_{\be-\rho,\vec j+\vec k}\ptl'',\be\ra.\eqno(3.23)$$
Since $\la\ptl',\be\ra=\la\ptl',\rho\ra\neq 0$ by the fact $\ptl'\in{\cal D}_{\be-\rho}\bs{\cal D}_{\rho}$, the above equation implies
$$\la e_{-\rho,\vec k}\ptl'',\be\ra=\la e_{\be-\rho,\vec j+\vec k}\ptl'',\be\ra\qquad\for\;\;\ptl''\in{\cal D}_{\rho}\cap{\cal D}_{\be}.\eqno(3.24)$$
Replacing $\be$ by $\be+\rho$ in (3.24), we obtain
$$\la e_{-\rho,\vec k}\ptl'',\be+\rho\ra=\la e_{\be,\vec j+\vec k}\ptl'',\be+\rho\ra\qquad\for\;\;\ptl''\in{\cal D}_{\rho}\cap{\cal D}_{\be}.\eqno(3.25)$$

On the other hand, taking the product (2.12) of (3.18) with $\be$ and choosing $\ptl''\in{\cal D}_{\be}\bs{\cal D}_{\rho}$, we have
$$\la e_{\be,\vec j}\ptl',\be\ra=0\qquad\for\;\,\ptl'\in{\cal D}_{\be}\cap{\cal D}_{\rho}.\eqno(3.26)$$
Taking $\be=\gm=-\rho$ and $\ptl',\ptl''\in{\cal D}_{\rho}$ in (3.10), we have
$$\la e_{-\rho,\vec j}\ptl',-\rho\ra\ptl''-\la e_{-\rho,\vec k}\ptl'',-\rho\ra\ptl'=0.\eqno(3.27)$$
By (2.14), $\dim{\cal D}_{\rho}\ge2$. Thus (3.27) forces
$$
\la e_{-\rho,\vec k}\ptl'',\rho\ra=0\qquad\for\,\;\ptl''\in{\cal D}_{\rho}.\eqno(3.28)$$
Substituting (3.26) and (3.28) into (3.25), we obtain
$$\la e_{-\rho,\vec k}\ptl'',\be\ra=\la e_{\be,\vec j+\vec k}\ptl'',\rho\ra\qquad\for\;\;\ptl''\in{\cal D}_{\rho}\cap{\cal D}_{\be}.\eqno(3.29)$$
Comparing this with (3.21), we have $\la e_{\be,\vec j}\ptl',\rho\ra=0$, which together with (3.26) implies that  the first term of (3.18) is zero. This completes the proof of (3.17).
\pse

{\it Claim 2}. There exists $\td{\ptl}$ such that
$$e_{\be,\vec j}\ptl'=\la\td{\ptl},\be\ra\ptl'-\la\ptl',\al\ra\td{\ptl}\qquad\for\;\,\ptl'\in{\cal D}_{\be-\rho},\;\be\in\G.\eqno(3.30)$$

Given an element  $\be\in\G\bs\BB{F}\al$, we can take $\hat{\ptl}_{\be}\in{\cal D}_{\be-\rho}$ such that
$\la\hat{\ptl}_{\be},\al\ra=1$. So we have a vector space
decomposition:
$${\cal D}_{\be-\rho}=\BB{F}\hat{\ptl}_{\be}\oplus({\cal D}_{\al}\cap{\cal D}_{\be}).\eqno(3.31)$$
Moreover, we define
$$\ptl_{\be,\vec j}=a_{\be,\vec j}\hat\ptl_{\be}-e_{\be,\vec j}\hat\ptl_{\be}.
\eqno(3.32)$$
Then for any $\ptl'=\ptl''+b\hat\ptl_{\be}\in{\cal D}_{\be-\rho}$ with $\ptl''\in{\cal D}_{\al}\cap{\cal D}_{\be}$, we have
$$e_{\be,\vec j}\ptl'=a_{\be,\vec j}\ptl''+b(a_{\be,\vec j}\hat\ptl_{\be}-\ptl_{\be,\vec j})=a_{\be,\vec j}\ptl'-\la\ptl',\al\ra\ptl_{\be,\vec j}\eqno(3.33)$$
for $\ptl'\in{\cal D}_{\be-\rho}$ and $\be\in\G\bs\BB{F}\al$ by (3.17). Using this in (3.10), we obtain
\begin{eqnarray*}& &(\la\ptl',\gm\ra(a_{\be}+a_{\gm}-a_{\be+\gm})+\la\ptl',\al\ra
(a_{\gm}-\la\ptl_{\be},\gm\ra))\ptl''-(\la\ptl'',\be\ra(a_{\be}+a_{\gm}-a_{\be+\gm})\\&&
+\la\ptl'',\al\ra(a_{\be}-\la\ptl_{\gm},\be\ra))\ptl'+\la\ptl',\al\ra\la\ptl'',\al\ra(\ptl_{\be}-\ptl_{\gm})+\la\ptl'',\be\ra\la\ptl',\al\ra(\ptl_{\be}-\ptl_{\be+\gm})\\& &+\la\ptl',\gm\ra\la\ptl'',\al\ra
(\ptl_{\be+\gm}-\ptl_{\gm})=0\hspace{8.8cm}(3.34)\end{eqnarray*}
for
$$\be,\gm,\be+\gm\in\G\bs\BB{F}\al,\ \ptl'\in{\cal D}_{\be-\rho},\;\ptl''\in{\cal D}_{\gm-\rho}.\eqno(3.35)$$

Let $\be,\gm\in\G$ such that $\{\al,\be,\gm\}$ is  linearly independent. First we take $\ptl'\in{\cal D}_{\al}\cap{\cal D}_{\be}\bs{\cal D}_{\gm}$ and $\ptl''\in{\cal D}_{\al}\cap{\cal D}_{\gm}\bs\{0\}$.
Then $\ptl'$ and $\ptl''$ are linearly independent. The last 3 terms in (3.34) vanish. So we obtain that the coefficient of $\ptl''$ is zero, which implies
$$a_{\be+\gm}=a_{\be}+a_{\gm}.\eqno(3.36)$$
By the above equation, we have
$$a_{\be}=a_{\be-\gm+\gm}=a_{\be-\gm}+a_{\gm}=a_{\be}+a_{-\gm}+a_{\gm},\eqno(3.37)$$
which implies
$$a_{-\gm}+a_{\gm}=0.\eqno(3.38)$$
Thus (3.36) and (3.38) imply
\begin{eqnarray*}\hspace{3cm}a_{n\be}&=&a_{n\be-\gm+\gm}\\ &=&a_{n\be-\gm}+a_{\gm}\\ &=&a_{\be}+a_{(n-1)\be-\gm}+a_{\gm}\\ &=& na_{\be}+a_{-\gm}+a_{\gm}\\ &=& na_{\be}\hspace{9.6cm}(3.39)\end{eqnarray*}
for $1\leq n\in\BB{Z}.$  Taking $\gm=\be$, $\ptl'\in{\cal D}_{\al}\cap{\cal D}_{\be}\bs\{0\}$ and
$\ptl''\in{\cal D}_{\be-\rho}\bs{\cal D}_{\al}$ in (3.34), we obtain
$$a_{\be}=\la\ptl_{\be},\be\ra.\eqno(3.40)$$
 Letting $\ptl'\in{\cal D}_{\be-\rho}\bs
{\cal D}_{\al}$ and $\ptl''\in{\cal D}_{\al}\cap{\cal D}_{\gm}\bs{\cal D}_{\be}$ in (3.34),
we get
$$(a_{\gm}-\la\ptl_{\be},\gm\ra)\ptl''+\la\ptl'',\be\ra(\ptl_{\be}-\ptl_{\be+\gm})=0
\eqno(3.41)$$
by (3.36).

Taking the product (2.12) of (3.41) with  $\al,\be$ and $\gm$, respectively, 
we obtain
$$\la\ptl_{\be},\al\ra=\la\ptl_{\be+\gm},\al\ra,\ \ a_{\gm}-\la\ptl_{\be},\gm\ra=
\la\ptl_{\be+\gm},\be\ra-\la\ptl_{\be},\be\ra,\ \ \la\ptl_{\be},\gm\ra=\la\ptl_{\be+\gm},\gm\ra.
\eqno(3.42)$$
By (3.40) and the last two equations of (3.42), we get
$$a_{\be}-\la\ptl_{\gm},\be\ra=\la\ptl_{\be+\gm},\gm\ra-\la\ptl_{\gm},\gm\ra=\la\ptl_{\be},\gm\ra-a_{\gm}.\eqno(3.43)$$
Replacing $\gm$ by $\gm-\be$ in the first equation of (3.42), we have
$$\la\ptl_{\be},\al\ra=\la\ptl_{\gm},\al\ra\;\for\;\be,\gm\in\G\;\mbox{such that}\;\{\al,\be,\gm\}\;\mbox{is linearly independent}.\eqno(3.44)$$
Taking the product (2.12) of (3.34) with  $\al$ and $\ptl',\ptl''\notin{\cal D}_{\al}$, we obtain
$$a_{\gm}-\la\ptl_{\be},\gm\ra=a_{\be}-\la\ptl_{\gm},\be\ra\eqno(3.45)$$
by (3.44). Comparing (3.45) with (3.43), we have that
the first term of (3.41) is zero. Thus replacing $\gm$ by $\gm-\be$ in (3.41), we get
$$\ptl_{\be}=\ptl_{\gm}\mbox{ for all $\be,\gm\in\G$ such that $\al,\be,\gm$ are linearly independent.}\eqno(3.46)$$

For any given $\be,\gm\in\G\bs\BB{F}\al$, we can choose
$\dlt\in\G$ such that $\{\al,\be,\dlt\}$ and $\{\al,\gm,\dlt\}$ are linearly independent. Thus we have 
$$\ptl_{\be}=\ptl_{\dlt}=\ptl_{\gm}.\eqno(3.47)$$
Hence $\ptl_{\be}$ does not depend on $\be\in\G\bs\BB{F}\al$, which is now denoted  by $\td{\ptl}$. Moreover, (3.33) and (3.40) imply
$$e_{\be,\vec j}\ptl'=\la\td\ptl,\be\ra\ptl'-\la\ptl',\al\ra\td\ptl
\ \ \for\,\;\be\in\G\bs\BB{F}\al,\ \ptl'\in{\cal D}_{\be-\rho}.\eqno(3.48)$$
For any $\gm=b\al\in\G$ and $\ptl'\in{\cal D}_{\gm-\rho}$, (3.10) and (3.48) show
\begin{eqnarray*}& &\la\be(\td{\ptl})\ptl'-\al(\ptl')\td\ptl,\gm\ra
\ptl''-\la\ptl'',\al+\be\ra(\be(\td\ptl)\ptl'-\al(\ptl')\td\ptl)\\& &+\la\ptl',\al+\gm\ra e_{\gm,\vec k}\ptl''-\la e_{\gm,\vec k}\ptl'',\be\ra\ptl'\\ &=& \la\td\ptl,\be+\gm\ra(\gm(\ptl')\ptl''-\be(\ptl'')\ptl')
-(\gm(\ptl')\al(\ptl'')-\be(\ptl'')\al(\ptl'))\td\ptl,\hspace{3cm}(3.49)\end{eqnarray*}
equivalently,
$$-\la e_{\gm,\vec k}\ptl''-\gm(\td\ptl)\ptl''+\al(\ptl'')\td\ptl,\be\ra\ptl'+(b+1)\al(\ptl')(e_{\gm,\vec k}\ptl''-\gm(\td\ptl)\ptl''+\al(\ptl'')\td\ptl)=0.\eqno(3.50)$$
Since $\ptl'\in{\cal D}_{\be-\rho}$ is arbitrary, we obtain the coefficient of $\ptl'$:
$$\la e_{\gm,\vec k}\ptl''-\gm(\td\ptl)\ptl''+\al(\ptl'')\td\ptl,\be\ra=0.\eqno(3.51)$$
By (2.12) and the arbitrariness of $\be\in\G\bs\BB{F}\al$, (3.51) implies
$$e_{\gm,\vec k}\ptl''-\gm(\td\ptl)\ptl''+\al(\ptl'')\td\ptl=0\ \ \for\,\ptl''\in{\cal D}_{\gm-\rho},\ \gm\in\G\cap\BB{F}\al. \eqno(3.52)$$
This together with (3.48) imply  (3.30).
\pse

Now we go back to (3.13). If $\al\ne\rho$, then $\td\ptl\in{\cal D}_{\al-\rho}$ and there exist $u'=\ol x^{\al,\vec i}\td\ptl\in {\cal S}_{\al}^{[\vec i]}$. It remains to consider the case $\al=\rho$.
\pse

{\it Claim 3}. If $\al=\rho$ and $i_{\ell_2}\ne0$, then $\chi_{\ell_2}(\td\ptl)=0$.
\pse

By assumption, $\rho=\al\ne0$ and $i_{\ell_2}\ne0$ imply that
$\ell_1=0$ and  $\ell_2>0$. By (3.30), $e_{\be,\vec j}$ does not depend on $\vec j$, which is simply denoted by $e_{\be}$. First we suppose that there exists $p\ne\ell_2$ such that $\rho_p\ne0$.
Take
$$u_1=D_{p,\ell_2}(t_{\ell_2})=t_{\ell_2}\ptl'+\ptl_p,\mbox{ where }\ptl'=\rho_p\ptl_{\ell_2}-\rho_{\ell_2}\ptl_p.\eqno(3.53)$$
We claim that the leading degree of $d(u_1)$ is $\le\vec i$. Otherwise, let $x^{\rho,\vec m}\ptl$ be the leading term of $d(u_1)$ with $\vec m>\vec i$ and $\ptl\in{\cal D}_0\bs\{0\}$. Take
$$u_2=x^{\rho}\ptl'',\mbox{ where }\ptl''=\ptl_p\mbox{ if }\rho(\ptl)=0\mbox{ or }\ptl''\in{\cal D}_0\cap{\cal D}_{\rho}\bs\{0\}\mbox{ otherwise}.\eqno(3.54)$$
Then $[d(u_1),u_2]$ has a term $x^{2\rho,\vec m}(\rho(\ptl)\ptl''-\rho(\ptl'')\ptl)\notin {\cal S}^{[\vec i]}$. On the other hand, $[u_1,u_2]\in {\cal S}^{[0]}$. So by (3.6), $d([u_1,u_2])\in {\cal S}^{[\vec i]}$ and $d(u_2)\in {\cal S}_{2\rho}^{[\vec i]}$. One can verify that $[u_1,{\cal S}_{2\rho}^{[\vec i]}]\subset {\cal S}_{2\rho}^{[\vec i]}$. This leads a contradiction. Hence we can assume
$$d(u_1)=x^{\rho,\vec i}\hat{\ptl}+\sum_{\vec n<\vec i}x^{\rho,\vec n}\ptl^{(\vec n)}\in {\cal S}_{\rho},
\eqno(3.55)$$
with $\hat{\ptl},\ptl^{(\vec n)}\in{\cal D}$. For any $\be\in\G$ with $\be(\ptl')=0$, we take $\ptl''\in{\cal D}_{\be-\rho}$. By (3.8) and considering the term with degree $\vec i$ in $d([u_1,x^{\be}\ptl''])$, we obtain
$$e_{\be}(\be_p\ptl''-\chi_{\ell_2}(\ptl'')\ptl')=\la\hat{\ptl},\be\ra\ptl''-\la\ptl'',\rho\ra\hat{\ptl}+(i_{\ell_2}\rho_p+\be_p)e_{\be}\ptl''-\chi_{\ell_2}(e_{\be}\ptl'')\ptl'.\eqno(3.56)$$
Applying (3.30) in (3.56) and noting that $\rho(\ptl')=0$, we have
\begin{eqnarray*}& &-\chi_{\ell_2}(\ptl'')\be(\td\ptl')\ptl'\\&=&\be(\hat\ptl)\ptl''-\rho(\ptl'')\hat\ptl
+(i_{\ell_2}+1)\rho_p(\be(\td\ptl)\ptl''-\rho(\ptl'')\td\ptl)\\& &-(\be(\td\ptl)\chi_{\ell_2}(\ptl'')-\rho(\ptl'')\chi_{\ell_2}(\td\ptl))\ptl',\hspace{8cm}(3.57)\end{eqnarray*}
equivalently,
$$((i_{\ell_2}+1)\rho_p\be(\td\ptl)+\be(\hat\ptl))\ptl''=\rho(\ptl'')(\hat\ptl+(i_{\ell_2}+1)\rho_p\td\ptl
-\chi_{\ell_2}(\td\ptl)\ptl')\;\;\: \for\,\be\in\G,\ \ptl''\in{\cal D}_{\be-\rho}.\eqno(3.58)$$
Since $\ptl''\in{\cal D}_{\be-\rho}$ is arbitrary, the coefficient of $\ptl''$ must be zero, which implies that the right-hand side in (3.58) is zero. We take $\be=\rho$ and $\ptl''=\ptl_p$ such that $\rho(\ptl'')\ne0$.
Applying $\chi_{\ell_2}$ to the right hand side of (3.58) and noting that $\chi_{\ell_2}(\ptl')=\rho_p$,
we obtain
$$\chi_{\ell_2}(\td\ptl)=-(i_{\ell_2}\rho_p)^{-1}\chi_{\ell_2}(\hat\ptl)=0,\eqno(3.59)$$
where the last equality follows from (3.55) and Lemma 2.2 and the assumption
that $i_{\ell_2}\ne0$.

Next we assume  that $\rho_p=0$ for all $p\in\ol{1,\ell}\bs\{\ell_2\}$. Then $\rho_{\ell_2}\ne0$. If $\ell_2\ge2$,  we take $u_1=D_{\ell_2,1}(t_1)=\rho_{\ell_2}t_1\ptl_1+\ptl_{\ell_2}$.  For any $\be\in\G$ with $\be_{\ell_2}=0$ and any $\ptl''\in{\cal D}_{\be-\rho}$,
we have
$$e_{\be}(\be_{\ell_2}\ptl''-\chi_1(\ptl'')\rho_{\ell_2}\ptl_1)=\la\hat\ptl,\be\ra\ptl''-\la\ptl'',\rho\ra\hat\ptl+(\rho_{\ell_2}+\be_{\ell_2})e_{\be}\ptl''-\chi_1(e_{\be}\ptl'')\rho_{\ell_2}\ptl_1
\eqno(3.60)$$
by the above arguments and considering $d([u_1,x^{\be}\ptl''])$.

By (3.30) and a similar argument as that for (3.57), the coefficient of $\ptl''$ must be zero. Then by calculating coefficient of $\ptl_{\ell_2}$, we obtain that $\chi_{\ell_2}(\td\ptl)=-\rho_{\ell_2}^{-1}\chi_{\ell_2}(\hat\ptl)=0$.

Finally, we assume that $\ell_2=1$ ($\ell_1=0$) and $\rho_p=0$ for all $p\in\ol{2,\ell}$. Then $\rho_1\ne0$, say $\rho_1=1$. For any $\ptl\in{\cal D}$, we have $\chi_1(\ptl)=\rho(\ptl)$. In this case, $\NJ=\BB{N}$ and $\vec j=j$ for all $\vec j\in\NJ$. For any $\be\in\G$, we define the linear map $\ol e_{\be}:{\cal D}_{\be-\rho}\rar{\cal D}$ by
$$d(x^{\be}\ptl')\equiv x^{\rho+\be,i}e_{\be}\ptl'+x^{\rho+\be,i-1}\ol e_{\be}\ptl'\ ({\rm mod\,}{\cal S}_{\rho+\be}^{(i-2)})\; \ \ \for\,\ptl'\in{\cal D}_{\be-\rho}\eqno(3.61)$$
(cf. (2.19)). By calculating the term of $d([x^{\be}\ptl', x^{\gm}\ptl''])$ with degree $i-1$, we obtain
\begin{eqnarray*}& &\la\ol e_{\be}\ptl',\gm\ra\ptl''-(i\chi_1(\ptl'')e_{\be}\ptl'+\la\ptl'',\rho+\be\ra
\ol e_{\be}\ptl')\\& &+(i\chi_1(\ptl')e_{\gm}\ptl''+\la\ptl',\rho+\gm\ra\ol e_{\gm}\ptl'')
-\la\ol e_{\gm}\ptl'',\be\ra\ptl'\\&=&\ol e_{\be+\gm}(\gm(\ptl')\ptl''-\be(\ptl'')
\ptl')\hspace{9.6cm}(3.62)\end{eqnarray*}
for $\be,\gm\in\G,\ \ptl'\in{\cal D}_{\be-\rho}$ and $\ptl''\in{\cal D}_{\gm-\rho}$ (cf. (3.10)).  Write
$\ol e_{\be}\ptl'$ as $\ol e_{\be}\ptl'=\hat{e}_{\be}\ptl'-i\be(\hat{\ptl})\ptl'$ for all $\be\in\G$ and $\ptl'\in{\cal D}_{\be-\rho}$. By (3.30) and the fact $\chi_1(\ptl)=\rho(\ptl)$, we observe that $\hat e_{\be}$ satisfies exactly the same equation as $e_{\be}$ in (3.10). The only difference between $\hat{e}_{\be}$ and $e_{\be}$ is that $\hat{e}_{\be}$ has codomain ${\cal D}$, while
$e_{\be}$ has codomain ${\cal D}_{\be}$ (cf. (3.7) with $\al=\rho$). However, this fact was not used when we used (3.10) in the proof of (3.30) if $\ell_1=0$. Therefore, we conclude that there exists $\hat\ptl\in{\cal D}$ such that
$$\hat e_{\be}\ptl'=\be(\hat\ptl)\ptl'-\rho(\ptl')\hat\ptl\qquad \for\;\,\ptl'\in{\cal D}_{\be-\rho},\ \be\in\G.\eqno(3.63)$$
Thus we can write (3.61) as
$$d(x^{\be}\ptl')\equiv x^{\rho+\be,i}(\be(\td\ptl)\ptl'-\rho(\ptl')\td\ptl)+x^{\rho+\be,i-1}(\be(\hat\ptl)\ptl'-\rho(\ptl')\hat\ptl-i\be(\td\ptl)\ptl')
\ ({\rm mod\,}{\cal S}_{\rho+\be}^{(i-2)}).\eqno(3.64)$$

Take $\be\in\G$ with $\be_1\ne0$. For any element
$$u=x^{\rho+\be,j}\ptl_u^{(j)}+x^{\rho+\be,j-1}\ptl_u^{(j-1)}+\cdots +x^{\rho+\be}\ptl_u^{(0)}
\in {\cal S}_{\rho+\be}^{[j]},\ \ \ j\ge0,\eqno(3.65)$$
we have
$$j\rho(\ptl_u^{(j)})+\be(\ptl_u^{(j-1)})=0.\eqno(3.66)$$
Note that any such element is a linear combination of $D_{p,q}(x^{\rho+\be,k})$ for $k\le j$ and $p,q\in\ol{1,\ell}$. Write $u=u_1+u_2$ such that $u_1$ is a linear combination of $D_{p,q}(x^{\rho+\be,j})$ and the leading degree of $u_2$ is $\le j-1$. Then by Lemma 2.2,
$\be(\ptl_{u_2}^{(j-1)})=0$, and we can assume
$$u_1=D_{p,q}(x^{\rho+\be,j})=x^{\rho+\be,j}(\be_q\ptl_p-\be_p\ptl_q)+jx^{\rho+\be,j-1}(\dlt_{q,1}\ptl_p-\dlt_{p,1}\ptl_q).\eqno(3.67)$$
By the fact $\rho_p=\dlt_{p,1},p\in\ol{1,\ell}$, we obtain (3.66).
Applying the result (3.66) to (3.64), we have
$$i(\be(\td\ptl)\rho(\ptl')-\rho(\ptl')\rho(\td\ptl))+(\be(\hat\ptl)\be(\ptl')-\rho(\ptl')\be(\hat\ptl)
-i\be(\td\ptl)\be(\ptl'))=0\; \ \for\,\ptl'\in{\cal D}_{\be-\rho}. \eqno(3.68)$$
By the fact $\be(\ptl')=\rho(\ptl')$ and taking $\ptl'\in {\cal D}_{\be-\rho}\bs{\cal D}_{\rho}$, we obtain $\rho(\td\ptl)=0$, that is, $\chi_{\ell_2}(\td\ptl)=0$. This proves Claim 3.

Next by our last claim, we can take some $u=\ol x^{\rho,\vec i}\td\ptl\in {\cal S}_{\rho}$ if $\vec i\ne\vec{0}$ and  $u=x^{\rho}\td\ptl\in {\cal W}_{\rho}^{[0]}$ if $\vec i=\vec{0}$. This completes the proof of Lemma 3.2.$\qquad\Box$
\psp

Now we shall describe homogeneous derivations of degree 0. Recall the derivations $\ptl_{t_i}$ and $\ptl^{\ast}_j$ on ${\cal A}={\cal A}(\ell_1,\ell_2,\ell_3;\G)$ defined in (1.11) for $i\in\ol{1,\ell_1+\ell_2}$ and $j\in\ol{1,\ell_2+\ell_3}$. Set
$${\cal D}^-=\sum_{i=1}^{\ell_2}\BB{F}\ptl_{t_{\ell_1+i}},\;\;{\cal D}^+=\sum_{j=1}^{\ell_2}\BB{F}\ptl^{\ast}_j,\;\;{\cal S}_0^+={\cal S}_0+{\cal D}+{\cal D}^+ +\sum_{p=1}^{\ell_1}\BB{F}t_p\ptl_p\eqno(3.69)$$
(cf. (1.11)). Note that ${\cal D}^-\subset {\cal S}^+_0$  by the second equation in (1.13).
Denote by ${\rm Hom}_{\BB{Z}}(\G,\BB{F})$ the set of additive group homomorphisms from $\G$ to $\BB{F}$.
\psp

{\bf Lemma 3.3}. {\it Let} $d\neq 0$ {\it be a derivation  satisfying  (3.6) with} $\al=0$.
{\it Then there exist} $u\in {\cal S}_0^+$ {\it and} $\mu\in\mbox{\it Hom}_{\BB{Z}}(\G,\BB{F})$ {\it 
such that} 
$$d(u')=[u,u']+\mu(\be)u'\qquad \for\;\,u'\in {\cal S}_{\be}.\eqno(3.70)$$

{\it Proof.} Obviously, (3.70) defines a homogeneous derivation of ${\cal S}$ of degree 0. Denote such
an derivation by $d_{u,\mu}$. We shall use the notations as before. In particular, $\vec i$ was defined in (3.7).  However, it is possible that $\vec i\notin\NJ$; for instance,
 $\vec i=(-1,0,\cdots,0)$ for $d=\ad_{\ptl_1}$.

For convenience, we shall assume $\ell_1+\ell_2\ge1$  because the case $\ell_1+\ell_2=0$ had been proved in [DZ] (cf. Theorem 4.2 there). Now (3.10) becomes
$$\la e_{\be,\vec j}\ptl',\gm\ra\ptl''-\be(\ptl'')e_{\be,\vec j}\ptl'+\gm(\ptl')e_{\gm,\vec k}-\la e_{\gm,\vec k}\ptl'',\be\ra\ptl'=e_{\be+\gm,\vec j+\vec k}(\gm(\ptl'')\ptl'-\be(\ptl'')\ptl'),
\eqno(3.71)$$
where $\be,\gm\in\G,\;\ptl'\in{\cal D}_{\be-\rho}$ and $\ptl''\in{\cal D}_{\gm-\rho}$
with the assumption (3.2).

Taking $\gm=0$ and $\vec k=0$ in (3.71),  we obtain $\la e_{0,0}\ptl'',\be\ra=0$ for all $\be\in\G$. Hence,
$$e_{0,0}\ptl''\in{\cal D}_1\qquad \for\;\,\ptl''\in{\cal D}_{\rho}.\eqno(3.72)$$
If $\ptl''\in{\cal D}_{\rho}\cap{\cal D}_3$, then for all
$u=\ol x^{\be,\vec j}\ptl'\in {\cal S}_{\be}$,  we have
$$\be(\ptl'')d(u)=d([\ptl'',u])=[d(\ptl''),u]+[\ptl'',d(u)]=[d(\ptl''),u]+\be(\ptl'')d(u).
\eqno(3.73)$$
due to that $\ad_{\ptl''}$ acts as the scalar $\be(\ptl'')$ on ${\cal S}_{\be}$.
This shows that $d(\ptl'')$ is in the center of ${\cal S}$. So
$$d(\ptl'')=0\qquad \for\;\,\ptl''\in{\cal D}_{\rho}\cap{\cal D}_3.\eqno(3.74)$$

For any $\ptl\in{\cal D}$, we define
$$p_{\ptl}=\left\{\begin{array}{ll}{\rm max}\{p\in\ol{1,\ell_1+\ell_2}\,|\,a_p\ne0\}&\mbox{if }\;\ptl=\sum_{p=1}^\ell a_p\ptl_p\in{\cal D}\bs{\cal D}_3,\\ \ell+1&\mbox{otherwise}.\end{array}\right.\eqno(3.75)$$
We shall choose a basis $B$ of ${\cal D}_{\rho}$ as follows. If $\rho\ne0$,
we take $p$ to be the largest index with $\rho_{p-\ell_1}\ne0$ when $\la{\cal D}_3,\rho\ra \ne0$ and the smallest index with $\rho_{p-\ell_1}\ne0$ when $\la{\cal D}_3,\rho\ra=0$ (cf. (2.12)). Let
$$B=\left\{\begin{array}{ll}\{\ptl_q\in{\cal D}_0\,|\,q\in\ol{1,\ell}\}&
\mbox{if }\rho=0,\\ \{\ptl_q-\rho_{p-\ell_1}^{-1}\rho_{q-\ell_1}\ptl_p\,|\,q\in\ol{1,\ell}\bs\{p\}\}&\mbox{if }\rho\ne0.\end{array}\right.\eqno(3.76)$$
By this choice of basis, we observe that different elements $\ptl\in
B\bs{\cal D}_3$ have different $p_{\ptl}$.
\pse

{\it Claim}. Replacing $d$ by some $d-d_{u,\mu}$ if necessary,
we can assume
$$e_{0,0}\ptl''=0\qquad \for\;\,\ptl''\in{\cal D}_{\rho}.\eqno(3.77)$$

Suppose $\ptl''\in B$ such that $e_{0,0}\ptl''\ne0$ with $p_{\ptl''}$
as minimal as possible. In this case, (3.72) and (3.74) force  $\ell_1\ge1,\ptl''\notin{\cal D}_3$
and $p=p_{\ptl''}\in\ol{1,\ell_1+\ell_2}$. We shall prove the claim by  induction on $p$.
For $\ptl'\in{\cal D}_{\be-\rho}$ and any $\be\in\G$ such that  $\be(\ptl'')=0$, we have
$$[\ptl'',\bar{x}^{\be,\vec j}\ptl']=j_p\bar{x}^{\be,\vec j-1_{[p]}}\ptl'.\eqno(3.78)$$
We want to prove that if $e_{0,0}\ptl''\ne0$, then we must have
$$i_q=0\qquad\for\;\;q\in\ol{p_{\ptl''}+1,\ell_1+\ell_2}\mbox{ \ and \ }e_{0,0}\ptl''\in{\cal D}^{(p)}={\rm Span}\{\ptl_1,...,\ptl_p\}\eqno(3.79)$$
(cf. (3.7) for $\vec i$). If $i_q\ne0$ for some $q>p$, then $\ell_1+\ell_2\ge2$, and so
$\ptl_q\in{\cal D}_{\rho-\rho}$ (cf. (2.13)). Taking $\be=\rho$ and $\ptl'=\ptl_q$ in (3.78), we have
$[\ptl'',x^{\rho}\ptl_q]=0$ and $[\ptl'',d(x^{\rho}\ptl_q)]\in {\cal S}_{\rho}^{[\vec i-1_{[p]}]}$.
However, the term with degree $\vec i-1_{[q]}$ in $[d(\ptl''),x^{\rho}\ptl_q]$ is $-i_qx^{\rho,\vec i-1_{[q]}}e_{0,0}\ptl''\notin {\cal S}^{[\vec i-1_{[p]}]}$, which leads a contradiction. Similarly, if $\chi_q(e_{0,0}\ptl'')\ne0$ for some $q\in\ol{p+1,\ell_1}$, then $[d(\ptl''),x^{\rho,1_{[q]}}\ptl_1]$ has a term
$\chi_q(e_{0,0}\ptl'')x^{\rho,\vec i}\ptl_1\in {\cal S}_{\rho}^{[\vec i]}\bs {\cal S}_{\rho}^{(\vec i)}$ and $[\ptl'',x^{\rho,1_{[q]}}\ptl_1]=0$. But the leading degree of $[\ptl'',d(x^{\rho,1_{[q]}}\ptl_1)]$
is $\le\vec i-1_{[p]}+1_{[q]}<\vec i$, which leads a contradiction again. This proves (3.79).

Since $e_{0,0}\ptl''\ne0$, we must have $\vec i\in\NJ$ because
otherwise ${\cal S}^{[\vec i]}={\cal S}^{(\vec i)}$, and we  have $e_{\be,0}\ptl''=0$ for all $\be\in\G$ by (3.8).
Thus $i_q\ge0$ for all $q\in\ol{1,\ell}$.

Let $p'=p_{ e_{0,0}\ptl''}$ (cf. (3.75)). By (3.72) and (3.79), $p'\le{\rm min}\{p,\ell_1\}$.
Write $e_{0,0}\ptl''=\sum_{q=1}^{p'}a_q\ptl_q$. According to the  assumption (2.34), $\rho=0$.
If $p'<\ell_1+\ell_2$, we let $k=i_{\ell_1+\ell_2}+1
+\dlt_{p,\ell_1+\ell_2}$ and take
\begin{eqnarray*}u_1&=&k^{-1}\sum_{q=1}^{p'}a_q D_{q,\ell_1+\ell_2}
(t^{\vec i+1_{[p]}+1_{[\ell_1+\ell_2]}})\\ &=&k^{-1}\sum_{q=1}^{p'}a_q(k t^{\vec i+1_{[p]}}
-(i_q+\dlt_{q,p})t^{\vec i+1_{[p]}+1_{[\ell_1+\ell_2]}-1_{[q]}})\\ &=&t^{\vec i+1_{[p]}}e_{0,0}\ptl''
-k^{-1}\sum_{q=1}^{p'}a_q(i_q+\dlt_{q,p})t^{\vec i+1_{[p]}+1_{[\ell_1+\ell_2]}-1_{[q]}}.\hspace{4.6cm}(3.80)\end{eqnarray*}

If $p'=\ell_1+\ell_2$, we must have $p'=p=\ell_1$ and $\ell_2=0$ due to $p'\le\ell_1$. Moreover,  (2.11) and Lemma 2.2 shows $\ell_1\ge2$ and $i_{\ell_1}=0$. By our choice of $B$ in (3.76), we obtain $\ptl''=\ptl_{\ell_1}$. According to the assumption of $p$, we get $e_{0,0}\ptl_q=0$ for $q\in\ol{1,\ell_1-1}$. We want to prove $\vec i=0$. Suppose  that $q'<\ell_1$ is the maximal index such that $i_{q'}\ne0$. Applying $d$ to 
$$(\dlt_{q,q'}-\dlt_{q,\ell_1})\ptl_q=[\ptl_q,t_{q'}\ptl_{q'}-t_{\ell_1}\ptl_{\ell_1}]
\ \ \ \for\,\;q\in\ol{q',\ell_1}
\eqno(3.81)$$
 and calculating the term with degree $\vec i$, we obtain
$$(\dlt_{q,q'}-\dlt_{q,\ell_1})e_{0,0}\ptl_q=\chi_{q'}(e_{0,0}\ptl_q)\ptl_{q'}-\chi_{\ell_1}(e_{0,0}\ptl_q)\ptl_{\ell_1}-i_{q'}e_{0,0}\ptl_q+(i_q+1)\ptl_{d(u_2)}^{(\vec i+1_{[q]})},
\eqno(3.82)$$
where $\ptl_{d(u_2)}^{(\vec i+1_{[q]})}$ is the element in ${\cal D}$ such that the term of $d(u_2)$ with degree $\vec i+1_{[q]}$ is $t^{\vec i+1_{[q]}}\ptl_{d(u_2)}^{(\vec i+1_{[q]})},$
and $u_2=t_{q'}\ptl_{q'}-t_{\ell_1}\ptl_{\ell_1}$.
Since $e_{0,0}\ptl_q=0$ for $q<\ell_1$, by (3.82), we have
$$\ptl_{d(u_2)}^{(\vec i+1_{[q]})}=0\qquad \for\;\,q\in\ol{q',\ell_1-1}.\eqno(3.83)$$
Since $i_q=0$ for $q\in\ol{q'+1,\ell_1}$, we get
$$\{\vec j\in\NJ\mid \vec i+1_{[q']}>\vec j\ge\vec i+1_{[\ell_1]}\}=\{ \vec i+1_{[q]}\mid q\in\ol{q'+1,\ell_1}\}.\eqno(3.84)$$
Thus by (3.83) and (3.84), we can assume
$$d(t_{q'}\ptl_{q'}-t_{\ell_1}\ptl_{\ell_1})\equiv t^{\vec i+1_{[\ell_1]}}\hat\ptl\;\;(\mbox{mod}\:{\cal S}^{[\vec{i}]})\eqno(3.85)$$
(cf. (2.20)), where $\hat\ptl=\ptl_{d(u_2)}^{(\vec i+1_{[\ell_1]})}\in{\cal D}_0$.
Now letting $q=\ell_1$ in (3.82) and applying $\chi_{\ell_1}$ to it (cf. (2.26)),
we have
$$-\chi_{\ell_1}(e_{0,0}\ptl_{\ell_1})=-\chi_{\ell_1}(e_{0,0}\ptl_{\ell_1})-i_{q'}\chi_{\ell_1}(e_{0,0}\ptl_{\ell_1})+\chi_{\ell_1}(\hat\ptl).
\eqno(3.86)$$
Thus by (3.85), (3.86) and Lemma 2.2, we obtain
$a_{\ell_1}=i_{q'}^{-1}\chi_{\ell_1}(\hat\ptl)=0$, which contradicts the case $p'=\ell_1$. Therefore $\vec i=0$. In this case, we take
$$u_1=t_{\ell_1}e_{0,0}\ptl''=a_{\ell_1}t_{\ell_1}\ptl_{\ell_1}+\sum_{q=1}^{\ell_1-1}a_qt_{\ell_1}\ptl_q\in {\cal S}_0^+.\eqno(3.87)$$

Now it is straightforward to verify that with the choices of $u_1$ in
(3.80) and (3.84), one has
$$u_1\in {\cal S}_0^+,\ [u_1,{\cal S}^{[\vec j]}]\subset {\cal S}^{[\vec i+\vec j]}\; \;\for\;\;
\vec j\in\NJ \;\mbox{ and }\;[u_1,\ptl]\subset {\cal S}^{(\vec i)}\eqno(3.88)$$
for $\ptl\in B$ with $\;p_{\ptl}<p$. Replacing $d$ by $d+(i_p+1)^{-1}\ad_{u_1}$, we have $e_{0,0}\ptl=0$ for all $\ptl\in B$ with $p_{\ptl}\le p$, which implies our claim (3.77) by induction on $p_{\ptl''}$.
\pse

Suppose $i_p<0$ for some $p\in\ol{1,\ell_1+\ell_2}$. Then for any $\vec j\in\NJ$ with $j_p<-i_p$, we have ${\cal S}^{[\vec i+\vec j]}={\cal S}^{(\vec i+\vec j)}$. Thus
$$e_{\be,\vec j}=0\qquad \for\;\,\be\in\G,\,\vec j\in\NJ \mbox{ with }j_p<-i_p\mbox{ if }i_p<0.\eqno(3.89)$$
Since $\ell\ge3$,   any element $\bar{x}^{\be,\vec k}\ptl'\in {\cal S}$ can be
generated by $\{\bar{x}^{\be,\vec j}\ptl'\mid\be\in\G,\,\vec j\in\NJ$ with $j_p<2\}$  for each $p\in\ol{1,\ell_1+\ell_2}$ by induction on $k_p$. This together with (3.89) shows that if $i_p<0$, then $i_p=-1$ and $i_q\ge0$ for any $q\in\ol{1,\ell_1+\ell_2}\bs\{p\}$.

For any $\ptl''\in{\cal D}_{\rho}$, we denote $p=p_{\ptl''}$. Let $\be\in\G$
with $\be(\ptl'')=0$ and let $\ptl'\in{\cal D}_{\be-\rho}\cap{\cal D}_{\rho}$.
Applying $d$ to (3.78), calculating the term with degree $\vec i+\vec j-1_{[p]}$ and using (3.77),
we obtain 
\begin{eqnarray*}j_p e_{\be,\vec j-1_{[p]}}\ptl'&=&\la\ptl_{d(\ptl'')}^{(\vec i-1_{[p]})},\be\ra\ptl'+\sum_{q=p+1}^{\ell_1+\ell_2}j_q\chi_q(\ptl_{d(\ptl'')}^{(\vec i-1_{[p]}+1_{[q]})})\ptl'
\\ & &-\sum_{q=p+1}^{\ell_1+\ell_2}\chi_q(\ptl')(i_q+1)\ptl_{d(\ptl'')}^{(\vec i-1_{[p]}+1_{[q]})}
+(i_p+j_p)e_{\be,\vec j}\ptl',\hspace{3.1cm}(3.90)\end{eqnarray*}
where $\ptl_{d(\ptl'')}^{(\vec i-1_{[p]})},
\ptl_{d(\ptl'')}^{(\vec i-1_{[p]}+1_{[q]})}$
are the elements in ${\cal D}$ such that the terms of $d(\ptl'')$ with degrees
$\vec i-1_{[p]}$, $\vec i-1_{[p]}+1_{[q]}$ are
$t^{\vec i-1_{[p]}}\ptl_{d(\ptl'')}^{(\vec i-1_{[p]})}$,
$t^{\vec i-1_{[p]}+1_{[q]}}\ptl_{d(\ptl'')}^{(\vec i-1_{[p]}+1_{[q]})}$, respectively.
Note that if $p=p_{\ptl''}=\ell+1$, this equation is trivial since all terms are zero.
Observe that on the right-hand side of (3.90), the second term does not
depend on $\be$ and the third terms does not depend on $\be,\vec j$. Moreover, $\ptl_{d(\ptl'')}^{(\vec i-1_{[p]})}=0$ if $i_p\le0$ or $i_{q'}<0$ for some $q'\ne p$ and
$\ptl_{d(\ptl'')}^{(\vec i-1_{[p]}+1_{[q]})}=0$ if $i_p\le0$ or $i_{q'}<0$ for some $q'\ne q$.

Using the arguments as those in the above and replacing $d$ by some $d-d_{u,\mu}$, we can assume
$$e_{0,\vec k}\ptl''=0\qquad \for\,\;\ptl''\in{\cal D}_{\rho},\;\vec k\in\NJ.\eqno(3.91)$$
Then setting $\be=0$ in (3.90) and using (3.91),
we see that the sum of the second term and the third term of the right-hand side is zero. Thus
$$j_p e_{\be,\vec j-1_{[p]}}\ptl'=\la\ptl_{d(\ptl'')}^{(\vec i-1_{[p]})},\be\ra\ptl'+(i_p+j_p)e_{\be,\vec j}\ptl'\eqno(3.92)$$
for any $\be\in\G,\,\ptl'\in{\cal D}_{\rho}\cap{\cal D}_{\be-\rho},\,
\ptl''\in{\cal D}_{\rho}\cap{\cal D}_{\be}$ with $p_{\ptl''}=p.$
By this, (3.71) and (3.91), we can show that if $\vec i\ne\vec{0}$, then
there exists $u=\bar{t}^{\vec i}\td\ptl\in {\cal S}_0$ such that
$$e_{\be,\vec j}\ptl'=\be(\td\ptl)\ptl'\qquad\for\;\;\be\in\G,\ \ptl'\in{\cal D}_{\be-\rho}.
\eqno(3.93)$$
Therefore the proof is completed by replacing $d$ by some $d-\ad_u$
and the induction on $\vec i$.

It remains to consider $\vec i=0$. In this case, by the proof of
Theorem 4.2 in [DZ], there exist an additive group homomorphism $\mu:\G\rar\BB{F}$ and
$\ptl\in{\cal D}_2+{\cal D}_3$ such that $d=d_{\ptl,\mu}.\qquad\Box$
\psp

{\bf Lemma 3.4}. {\it Every homogeneous derivation} $d\in (\der{\cal S})_{\al}$ {\it 
must satisfy the condition in (3.6)}.
\psp

{\it Proof}. Choose $\G'\subset\G$ to be a nondegenerate subgroup of $\G$ generated by a
finite subset $\G'_0$ of $\G$ such that $\rho,\al\in\G'_0$.
Let ${\cal S}'$ be the Lie subalgebra of ${\cal S}$ generated by
$$\{D_{p,q}(x^{\be,\vec j})\,|\,\be\in\G'_0,\vec j\in\NJ,\,|\vec j|\le4,\,p,q\in\ol{1,\ell}\}.\eqno(3.94)$$
Then it is straightforward to check that
${\cal S}'={\cal S}(\ell_1,\ell_2,\ell_3;\rho,\G')$,
and $d'=d|_{{\cal S}'}$ is a homogeneous derivation of ${\cal S}'$ of degree $\al$.
Since (3.94) is a finite set and a  derivation is determined by its action on generators, we see that  the derivation $d'$ of ${\cal S}'$ satisfies the condition in (3.6).

For $u\in{\cal S}$ and $\mu\in\mbox{Hom}_{\BB{Z}}(\G,\BB{F})$, we use the notation $d_{u,\mu}$ to denote the derivation defined by the right-hand side of the equation (3.70). By Lemmas 3.1-3, there exist $u'=u'_1+u'_2\in {\cal S}'_{\al}+{\cal W}^{[0]}_{\rho}$ and $\mu'\in{\rm Hom}_{\BB{Z}}(\G',\BB{F})$ such that $d'=d_{u',\mu'}$, $u'_2=0$ if $\al\ne\rho$ and $\mu'=0$ if $\al\ne0$.

We claim that for $u\in{\cal S}$ and $\mu\in{\rm Hom}_{\BB{Z}}(\G,\BB{F})$,
$d_{u,\mu}|_{{\cal S}_{\G_1}}=0$ if and only if $u\in{\cal D}_3$ and
$\mu(\al)=-\al(u)$ for any $\al\in\G_1$, where
$\G_1$ is any nondegenerate subgroup of $\G$ and
${\cal S}_{\G_1}={\cal S}(\ell_1,\ell_2,\ell_3;\rho,\G_1)$. Write
$$u=\sum_{(\al,\vec i)\in K}x^{\al,\vec i}\ptl^{(\al,\vec i)},\qquad\mbox{where}\;\;\ptl^{(\al,\vec i)}\in{\cal D}\;\;\mbox{and}\eqno(3.95)$$
$$K=\{(\al,\vec i)\in \G\times\NJ\,|\,\ptl^{(\al,\vec i)}\ne0\}\eqno(3.96)$$
 is a finite set. If some $(\al,\vec i)\in K$ with $\vec i\ne \vec{0}$ or $\al\notin\BB{F}\rho$,
then clearly, there exists some $\ptl\in{\cal D}_{\rho}$ such that
$$d_{u,\mu}=[u,\ptl]=-\sum_{(\al,\vec i)\in K}\ptl(x^{\al,\vec i})\ptl^{(\al,\vec i)}\ne0.\eqno(3.97)$$ Thus we can rewrite
$$u=\sum_{c\in K'}x^{c\rho}\ptl^{(c\rho)},\qquad\mbox{where}\;\;K'=\{c\in\BB{F}\,|\,c\rho\in\G,\ptl^{(c\rho)}\ne0\}\;\;\mbox{is a finite set}.\eqno(3.98)$$
If $\chi_p(\ptl^{(c\rho)})\ne0$ for some $p\in\ol{1,\ell_1+\ell_2}$,
then we can choose $\ptl\in{\cal D}_{0}\bs\{0\}$ with $\chi_p(\ptl)=0$. So $t_p\ptl\in {\cal S}$ and 
$$d_{u,\mu}(t_p\ptl)=[u,t_p\ptl]\ne0.\eqno(3.99)$$
Thus all $\ptl^{(c\rho)}\in{\cal D}_3$ (cf. (2.1)). Let $0\neq c\in K'$. Take $\be\in\G_1$ with $\be(\ptl^{(c\rho)})\ne0$ and choose $\ptl\in
{\cal D}_{\be}\cap{\cal D}_{\rho}$. Then we have
$$d_{u,\mu}(x^{\be}\ptl)=\sum_{c\in K'}\be(\ptl^{(c\rho)})x^{c\rho+\be}\ptl^{(c\rho)}
+\mu(\be)x^{\be}\ptl\ne0.\eqno(3.100)$$
Hence $u=\ptl^{(0)}\in{\cal D}_3$. Then for any $x^{\al}\ptl\in {\cal S}_{\G_1}$, we have
$$0=d_{u,\mu}(x^{\al}\ptl)=(\al(u)+\mu(\al))x^{\al}\ptl,\eqno(3.101)$$
which implies $\mu(\al)=-\al(u)$ for all $\al\in\G_1$. This proves the claim.

We can regard ${\cal D}_3$ as a subspace of ${\rm Hom}_{\BB{Z}}(\G,\BB{F})$. Choose a subspace
${\rm Hom}^*_{\BB{Z}}(\G,\BB{F})$ of ${\rm Hom}_{\BB{Z}}(\G,\BB{F})$ such that
${\rm Hom}_{\BB{Z}}(\G,\BB{F})={\cal D}_3\oplus{\rm Hom}^*_{\BB{Z}}(\G,\BB{F})$ as vector
spaces. Thus we can always assume $\mu\in{\rm Hom}^*_{\BB{Z}}(\G,\BB{F})$ when we use the notation $d_{u,\mu}$. Hence for any $u\in {\cal S}_{\al}+{\cal W}_{\rho}^{[0]}$ and  $\mu\in{\rm Hom}^*_{\BB{Z}}(\G,\BB{F})$, $d_{u,\mu}|_{{\cal S}_{\G_1}}=0$ implies $u=0$ and $\mu|_{\G_1}=0$.

Let $\G_2$ be the maximal subgroup of $\G$ such that for ${\cal S}_{\G_2}={\cal S}(\ell_1,\ell_2,\ell_3;\G_2)$, there exists $\mu_2\in{\rm Hom}_{\BB{Z}}(\G_2,\BB{F})$ with $d|_{{\cal S}_{\G_2}}=d_{u_1,\mu_2}$ and $\mu_2|_{\G_1}=\mu_1$.

Suppose $\G_2\ne\G$. Take $\be\in\G\bs\G_2$. Let $\G_3$ be the subgroup of $\G$ generated by $\G'$ and $\be$. Then $\G_3$ is finitely generated. Thus there exist $u_3,\mu_3$ such that
$d|_{{\cal S}_{\G_3}}=d_{u_3,\mu_3}$. Then we have
$d_{u_3-u_1,\mu_3-\mu_1}|_{{\cal S}_{\G'}}=d|_{{\cal S}_{\G'}}-d|_{{\cal S}_{\G'}}=0$. Thus
$u_3=u_1$ and $\mu_3|_{\G'}=\mu_1|_{\G'}$. Similarly, $\mu_2|_{\G_2\cap\G_3}=\mu_3|_{\G_2\cap\G_3}$.

Let $\G_4$ be the subgroup of $\G$ generated by $\G_2$ and $\be$. Define $\mu\in{\rm Hom}_{\BB{Z}}(\G_4,\BB{F})$ as follows. For any $\gm\in\G_4$, we can write $\gm=\tau+n\be$ with $n\in\BB{Z}$ and $\tau\in\G_2$. Define $\mu_4(\gm)=\mu_2(\tau)+n\mu_3(\be)$. Suppose $\tau+n\be=0$ for some $n\in\BB{Z}$. Then $\tau=-n\be\in\G_2\cap\G_3$. Since $\mu_2|_{\G_2\cap\G_3}=\mu_3|_{\G_2\cap\G_3}$, we have $\mu_2(\tau)=\mu_3(\tau)=\mu_3(-n\be)$.
But obviously, $\mu_3(-n\be)=-n\mu_3(\be)$. Hence $\mu_2(\tau)+n\mu_3(\be)=0$. This shows that
$\mu_4\in{\rm Hom}_{\BB{Z}}(\G_4,\BB{F})$ is uniquely defined. So $d|_{{\cal S}_{\G_4}}=d_{u_1,\mu_4}$ and $\G_4\supset\G_2$, $\G_4\ne\G_2$. This contradicts the maximality of $\G_2$. Therefore, $\G_2=\G$ and $d=d_{u_1,\mu_2}$ satisfies the condition in (3.6).$\qquad\Box$
\psp

{\bf Lemma 3.5}. {\it Let} $d$ {\it be any derivation of} ${\cal S}$. {\it Write}
$$d=\sum_{\al\in \G}d_{\al}\qquad\mbox{\it with}\;\; d_{\al}\in(\mbox{\it Der}\:{\cal S})_{\al}.\eqno(3.102)$$
{\it Then}
$$d_{\al}=0\;\;\mbox{\it for all but a finite}\;\;\al\in\G.\eqno(3.103)$$

{\it Proof}. By Lemmas 3.1-4, for any $\al\in\G\bs\{0,\rho\}$, there exists
$u_{\al}=\bar{x}^{\al,\vec i_{\al}}\ptl_{\al}\in {\cal S}$ such that $d_{\al}=\ad_{u_{\al}}$.
We always assume that $\ptl_{\al}\ne0$ if $d_{\al}\ne0$ and $\vec i_{\al}$ is the leading degree of $u_{\al}$. We shall prove that
$$Y=\{\al\in\G\bs\{0,\rho\}\mid\ptl_{\al}\ne0\}\eqno(3.104)$$
is finite. Take a $\BB{F}$-basis $\{\al^{(\ell_1+1)},...,\al^{(\ell)}\}$ of $\BB{F}^{\ell_2+\ell_3}$ from $\G\bs\{\rho\}$. Define
$$Z'=\{\al\in Y\mid\ptl_{\al}\in{\cal D}_1\},\;\;Z_{\be}=\{\al\in Y\mid\be(\ptl_{\al})\ne0\}
\ \ \for\,\be\in\G.\eqno(3.105)$$
Then
$$Y=Z'\bigcup\bigcup_{p\in\ol{\ell_1+1,\ell}}Z_{\al^{(p)}}.\eqno(3.106)$$
 So it is sufficient to prove that $Z'$ is a finite set and so is $Z_{\be}$ for any $\be\in\G\bs\{0,\rho\}$.

Suppose that $Z'$ is an infinite set. Then $\ell_1>0$, and so we can assume $\rho=0$ by Lemma 2.3. Hence ${\cal D}_0\supset{\cal D}_2+{\cal D}_3$, where $0\in\G$. Take $\ptl\in{\cal D}_2+{\cal D}_3$ such that there are infinite many $\al\in Z'$ with $\al(\ptl)\ne0$. For any $\al\in Z'$, we have
$$d_{\al}(\ptl)\equiv-\al(\ptl)\bar{x}^{\al,\vec i_{\al}}\ptl
\ ({\rm mod\ }{\cal S}_{\al}^{(\vec i_{\al})}).\eqno(3.107)$$
Thus there are infinitely many $\al$ with $d_{\al}(\ptl)\ne0$, which contradicts the fact that $d(\ptl)$ is contained in a sum of finite number of ${\cal S}_{\al}$.

Assume that $\al\in Z_{\be}$. For $\ptl\in{\cal D}_{\be-\rho}$,
we have
$$d_{\al}(\ol x^{\be,\vec j}\ptl)\equiv
\bar{x}^{\al+\be,\vec i_{\al}+\vec j}(\be(\ptl_{\al})\ptl-\al(\ptl)\ptl_{\al})
\ ({\rm mod\ }{\cal S}_{\al+\be}^{(\vec i_{\al}+\vec j)})\eqno(3.108)$$
(cf. (2.19)). The fact that $d(\ptl)$ is contained in a sum of finite number of ${\cal S}_{\al}$ implies that
$$\be(\ptl_{\al})\ptl-\al(\ptl)\ptl_{\al}=0\eqno(3.109)$$                                 
for all but a finite $\al\in Z_{\be}$. Since $\be(\ptl_{\al})\ne0$ by (3.105), $\ptl\in{\cal D}_{\be-\rho}$ is arbitrary and $\dim{\cal D}_{\be-\rho}\ge2$, (3.109) implies that $\ptl_{\al}=0$ for all but a finite $\al\in Z_{\be}$. Therefore, $Z_{\be}$ is a finite set.$\qquad\Box$
\psp

By Lemma 3.5, we have
$$\der {\cal S}=\bigoplus_{\al\in\G}(\der{\cal S})_{\al}
\eqno(3.110)$$
(cf. (3.5)). For convenience, we identify an additive function $\mu$ with the derivation
$d_{0,\mu}$ defined by $\mu$. Thus
$${\rm Hom}_{\BB{Z}}(\G,\BB{F})\subset (\der{\cal S})_0.\eqno(3.111)$$
Recall the notations in (3.69). For any $\ptl\in{\cal D}_2+{\cal D}_3$, there is a unique way
to decompose
$$\ptl=\ptl^+ +\ptl^-\qquad\mbox{with}\;\;\ptl^+\in{\cal D}^+ +{\cal D}_3,\;\ptl^-\in{\cal D}^-.
\eqno(3.112)$$
For any $\al\in\G$, we define
$$\la\ptl^+,\al\ra=\la\ptl,\al\ra
\eqno(3.113)$$
(cf. (2.12)). Then ${\cal D}^+ +{\cal D}_3$ can be identified with a subspace of $\BB{F}$-linear function ${\rm Hom}_{\BB{F}}(\G,\BB{F})$. Thus
$${\cal D}^+ +{\cal D}_3\subset
{\rm Hom}_{\BB{Z}}(\G,\BB{F})\subset (\der{\cal S})_0.\eqno(3.114)$$
We shall identify $u$ with $\ad_u|_{\cal S}$ for $u\in {\cal W}$ when the context is clear. In
particular, for $\mu\in{\rm Hom}_{\BB{Z}}(\G,\BB{F})$ and $u\in {\cal S}_{\al}$, we have
$$ [\mu,u]=\mu(\al)u.\eqno(3.115)$$
We summarize the results in Lemmas 3.1-3.5 as the following theorem.
\psp

{\bf Theorem 3.6}. {\it The derivation algebra} $\der{\cal S}$ {\it is an} $\G$-{\it graded Lie algebra (cf. (3.5) and (3.110)) with}
$$(\der{\cal S})_{\al}=\left\{\begin{array}{ll}{\cal S}_{\al}&\mbox{\it if }\al\ne\rho,0,\\ {\cal S}_{\rho}+{\cal W}_{\rho}^{[0]}&\mbox{\it if }\al=\rho\ne0,
\\ {\cal S}_0+\BB{F} t_{\ell_1}\ptl_{\ell_1}+{\cal D}+{\rm Hom}_{\BB{Z}}(\G,\BB{F})&\mbox{\it if }\al=0.\end{array}\right.
\eqno(3.116)$$

\section{Proof of the Main Theorem}

In this section, we shall present the proof of the main theorem in this paper. First we need three more lemmas on derivations.

A linear transformation $T$ on a vector space $V$ is called {\it locally-nilpotent} if for any $v\in V$, there exists a positive integer $n$ (depending on $v$) such that
$$T^n(v)=0.\eqno(4.1)$$
\psp

{\bf Lemma 4.1}. {\it If} $d\in\der{\cal S}$ {\it is locally-finite (cf. (1.16)), then}
$$
d\in {\cal A}{\cal D}_1+{\cal D}+{\rm Hom}_{\BB{Z}}(\G,\BB{F}).\eqno(4.2)$$
{\it If} $d$ {\it is locally-nilpotent, then}
$$d\in {\cal A}{\cal D}_1\oplus{\cal D}^-.\eqno(4.3)$$

{\it Proof}. If $\ell_2+\ell_3=0$, there is nothing to prove. Suppose $\ell_2+\ell_3\ge1$.
Choose a total order $\le$ on $\G$ compatible with its group structure.
Let $d$ be a locally-finite derivation. By Lemma 3.5, there exists a finite subset $\G_0$ of $\G$ such that
$d=\sum_{\al\in\G_0}d_{\al}$. Write $d$ as
$$d=u+\ptl^-+\mu,\;u=\sum_{(\al,\vec i)\in\G_0\times\vec J_0}x^{\al,\vec i}\ptl^{(\al,\vec i)},\eqno(4.4)$$
where $\G_0\times\vec J_0$ is a finite subset of $\G\times\NJ$ and $\ptl^{(\al,\vec i)}\in{\cal D},$ $\ptl^-\in{\cal D}^-$, $\mu\in{\rm Hom}_{\BB{Z}}(\G,\BB{F}).$ We want to prove
$$\left\{\begin{array}{ll}\ptl^{(\al,\vec i)}\in{\cal D}_1\;\;\for\,(0,0)\ne(\al,\vec i)\in\G_0\times\vec J_0 &\mbox{if $d$ is locally-finite;}\\ \ptl^{(\al,\vec i)}\in{\cal D}_1\;\;\for\,(\al,\vec i)\in\G_0\times\vec J_0\mbox{ and }\mu=0 &\mbox{if $d$ is locally-nilpotent}.\end{array}\right.\eqno(4.5)$$
We shall only give the proof of the first statement in (4.5). The proof of the second statement in (4.5) is similar. Suppose that $\ptl^{(\dlt,\vec i)}\notin{\cal D}_1$ for some $\dlt\in\G_0\bs\{0\}$
(the case $\vec i\ne\vec{0}$ can be similarly proved). Take
$$\dlt={\rm max}\{\al\in\G_0\,|\,\ptl^{(\al,\vec i)}\notin{\cal D}_1\mbox{ for some }\vec i\mbox{ with }(\al,\vec i)\in\G_0\times\vec J_0\}\eqno(4.6)$$
and
$$\vec j={\rm max}\{\vec i\in\vec J_0\,|\,(\dlt,\vec i)\in\G_0\times\vec J_0\mbox{ with }\ptl^{(\dlt,\vec i)}\notin{\cal D}_1\}.\eqno(4.7)$$
Note that for any $\al\in\G$ and $p\in\ol{1,\ell_1}$, $x^{\al}\ptl_p$ is a locally-nilpotent derivation on ${\cal S}$. Hence ${\rm exp}(x^{\al}\ptl_p)$ is an automorphism of ${\cal S}$. Let $G$ be the subgroup of ${\rm Aut}({\cal S})$ generated by such automorphisms. Elements in $G$ induce automorphisms of $\der{\cal S}$.
\pse

{\it Claim 1}. Replacing $d$ by $g(d)$ (which is again
locally-finite) for some $g\in G$, we can assume 
$$j_p=0\;\;\; \for\,\;p\in\ol{1,\ell_1},\eqno(4.8)$$
where $\vec j$ is defined in (4.7).
\pse

We define another total order $>'$ on $\NJ$ different from (2.16) by
$$\vec i >'\vec j\ \lra\mbox{ for the first $p$ with $i_p\ne j_p$, we have $i_p>j_p$}.
\eqno(4.9)$$
With respect to this order, we define 
$$\vec k={\rm max}\{\vec i\in\vec J_0\,|\,\ptl^{(\al,\vec i)}\notin{\cal D}_1\mbox{ for some }\al\in\G_0\mbox{ with }(\al,\vec i)\in\G_0\times\vec J_0\},\eqno(4.10)$$
$$\tau={\rm max}\{\al\in\G_0\,|\,(\al,\vec k)\in\G_0\times\vec J_0\mbox{ with }\ptl^{(\al,\vec k)}
\notin{\cal D}_1\}.\eqno(4.11)$$
By this definition, we have
$$\ptl^{(\al,\vec i)}\notin{\cal D}_1\lra(\al,\vec i)=(\tau,\vec k)\mbox{ or }
\vec i<'\vec k\mbox{ or }\vec i=\vec k\mbox{ but }\al<\vec\tau.\eqno(4.12)$$
For each $p\in\ol{1,\ell_1}$, we take
$$m_p=\sum_{q=p+1}^{\ell_1}\ \sum_{(\al,\vec i)\in\G_0\times\vec J_0}i_q,
\eqno(4.13)$$
and $\al^{(p)}\in\G$ such that
$$\al^{(p)}>m_p\al^{(p+1)}\;\;\for\;\; p\in\ol{1,\ell_1-1}\mbox{ and }\al^{(\ell_1)}>\be,
\eqno(4.14)$$
where $\be$ is the largest element of $\G_0$.
Set
$$d'=\prod_{p=1}^{\ell_1}{\rm exp}(x^{\al^{(p)}}\ptl_p)(d).\eqno(4.15)$$
Then  $d'$ has a term 
$$x^{\mu,\vec k'}\ptl^{(\tau,\vec k)}\mbox{ with }\mu=\sum_{p=1}^{\ell_1}k_p\al^{(p)}+\tau,\
\vec k'=\vec k-\sum_{p=1}^{\ell_1}(k_p)_{[p]},\;\;\ptl^{(\vec\tau,\vec k)}\notin{\cal D}_1,
\eqno(4.16)$$
where $\tau$ and $\vec k$ are defined in (4.10) and (4.11), and any term
$x^{\al',\vec i'}\ptl^{(\al',\vec i')}$ appears in $d'$ with $\ptl^{(\al',\vec i')}\notin{\cal D}_1$ must be of the form
$x^{\nu,\vec m}\ptl^{(\al,\vec i)}$ with
$$\nu=\sum_{p=1}^{\ell_1}n_p\al^{(p)}+\al,\;n_p\le i_p\;\,\for\;\,p\in\ol{1,\ell_1}\mbox{ and }
\vec m=\vec i-\sum_{p=1}^{\ell_1}(n_p)_{[p]},\
\ptl^{(\al,\vec i)}\notin{\cal D}_1.\eqno(4.17)$$

Let $\nu$ and $\vec i$ be as in (4.17). By (4.10) and (4.11), we have the following two cases. First we have $\vec i<'\vec k$. We let $p'$ be the first index with $i_{p'}<k_{p'}$. By  (4.12) and (4.13),
$$\nu\leq\sum_{p=1}^{p'}n_p\al^{(p)}+m_{p'}\al^{(p'+1)}<\sum_{p=1}^{p'-1}k_p\al^{(p)}+k_{p'}\al^{(p')}
\leq\mu.\eqno(4.18)$$
Secondly we have $\vec i=\vec k$ but $\al<\tau$, which implies $\nu<\mu$. Thus we have proved that $\mu$ in (4.16) is
the element $\dlt$ defined in (4.6) for $d'$ and $\vec k'$ is the corresponding element $\vec j$ defined in (4.7) for $d'$, which satisfies (4.8) by (4.16). This proves Claim 1.
\pse

{\it Claim 2}. Replacing $d$ by $\ol\psi(d)$ for some
$\ol\psi$ defined in (2.31) and (2.32), and considering the derivation
$\ol\psi(d)$ of ${\cal S}'$, we can assume that if a term $x^{\al,\vec i}\ptl_p$
appears in $d$ with $\al\geq\dlt$ and $p\in\ol{1,\ell_1}$, then
$$p>1\mbox{ and }i_1=0,\mbox{ or }p=1\mbox{ and } i_1\le1,i_q=0\,\;\for\ q\in\ol{2,\ell_1}.
\eqno(4.19)$$
\pse

Pick $\mu,\nu\in\G$ such that
$$\mu>(m_1+1)\nu,\ \ \nu>\be\mbox{ and } \nu>-{\rm min}\{\al\,|\,(\al,\vec i)\in\G_0\times\vec J_0\},
\eqno(4.20)$$
where $m_1$ is defined in (4.13) with $p=1$. Define $\ol\psi$ as in (2.31) and (2.32) with
$$\al^{(1)}=-\mu,\ \al^{(p)}=-\nu\;\; \for\;\;p\in\ol{2,\ell_1}.\eqno(4.21)$$
Then we have
$$\ol\psi(x^{\al,\vec i}\ptl_p)=x^{\gm,\vec i}\ptl_p\;\;\mbox{with}\;\;\gm=\al-i_1\mu-\sum_{p=2}^{\ell_1}i_p\nu+
\dlt_{p,1}(\mu-\nu)+\nu\;\;\for\;\;p\in\ol{1,\ell_1},\eqno(4.22)$$
$$\ol\psi(x^{\al,\vec i}\ptl_p)\in {\cal S}_{\gm}\;\;\mbox{with}\;\;\gm=\al-i_1\mu-\sum_{p=2}^{\ell_1}i_p\nu\;\;\for\;\;p\in\ol{\ell_1+1,\ell}.\eqno(4.23)$$
By (4.22) and (4.23),  we see that any term $x^{\al,\vec i}\ptl_p$ appearing in $d$ is mapped to ${\cal S}_{\gm}$
with some $\gm<\dlt$ if $\vec i$ does not satisfy (4.19). On the other hand, if we write $\ptl^{(\dlt,\vec j)}=\ptl'+\ptl''$ with $\ptl'\in{\cal D}_1$, $\ptl''\in({\cal D}_2+{\cal D}_3)\bs\{0\}$ for the term $x^{\dlt,\vec j}\ptl^{(\dlt,\vec j)}$ in $d$, then by (2.32), $x^{\dlt,\vec j}\ptl''$ is mapped to an element in
${\cal S}_{\dlt}$ with a term $x^{\dlt,\vec j}\ptl''$. Thus we have (4.19) for $\ol\psi(d)$.
\pse

Now suppose $d$ satisfies (4.8) and (4.19) with $\dlt,\vec j$ defined by (4.6) and (4.7).
If there exists a term $x^{\al,\vec i}\ptl_1$ appearing in $d$ with $\al>\dlt$ and $i_1=1$,
then we take $\dlt'$ to be the largest such $\al$, $\vec j$ to be such $\vec i$ and
take $v=x^{\gm}\ptl_1$ with $\gm=0$. Otherwise we take $\dlt'=\dlt$ and $v=x^{\gm}\ptl_1$ with $\gm\in\G$ such that $\la\ptl^{(\dlt,\vec j)},\gm\ra\ne0$. Then by (4.19), it is straightforward to verify that
$$d^k(v)\equiv ad_u^k(v)\not\equiv 0\,({\rm mod}\;{\cal S}_{(\gm+k\dlt')})\;\;\mbox{with}\;\;u=x^{\dlt',\vec j}\ptl^{(\dlt',\vec j)},\eqno(4.24)$$
where ${\cal S}_{(\al)}=\bigoplus_{\al>\gm\in\G}{\cal S}_\gm$.
Thus $d$ is not locally-finite. This proves the first statement in (4.5), and the second can be proved similarly.$\qquad\Box$
\psp

An element $u\in {\cal W}$ is called {\it locally-nilpotent} if $\ad_u$ is a locally-nilpotent. If $\ell_1\ge1$ in the proof of Lemma 4.1, then we can always take locally-nilpotent $v\in {\cal A}{\cal D}_1$ in (4.24). Denote
$${\rm Nil}({\cal S})=\mbox{ the subalgebra  generated by locally-nilpotent
elements in ${\cal S}\cap {\cal A}{\cal D}_1$}.\eqno(4.25)$$
Then we have the following lemma.
\psp

{\bf Lemma 4.2}. {\it Suppose} $\ell_1\ge1$.
{\it If} $d\in\der {\cal S}$ {\it is locally-nilpotent on} ${\rm Nil}({\cal S})$, {\it then} $d\in {\cal A}{\cal D}_1\oplus{\cal D}^-$.$\qquad\Box$
\psp

Set
$$A_1={\rm Span}\{x^{\al}\,|\,\al\in\G\},\;\; A_0=\BB{F}[t_1,t_2,...,t_{\ell_1+\ell_2}].\eqno(4.26)$$
We have ${\cal A}=A_0\cdot A_1$.
\pse

{\bf Lemma 4.3}. {\it Let} $u=x^{\al}\in A_1,\,
d=\bar{x}^{\be,\vec j}\ptl\in {\cal S}_{\be}\cap {\cal A}{\cal D}_1$ {\it be homogeneous elements}.

{\it (i) If} $\ad_d$ {\it is locally-nilpotent on} ${\cal S}\cap {\cal A}{\cal D}_1$, {\it then} $\ad_{ud}$ {\it is locally-nilpotent on} ${\cal A}{\cal D}_1$.

{\it (ii) If} $\al+\be\ne\rho$ {\it or} $\ell_2\ge1$, {\it  then} $ud\in {\cal S}$.
\psp

{\it Proof}. (i) Note that for any $v=x^{\be,\vec i}\ptl\in {\cal A}{\cal D}_1$, we have
$$(\ad_{ud})^m(v)=u^m(\ad_d)^m(v)\ \ \for\,m\in\BB{N}.\eqno(4.27)$$
This proves (i).

(ii) Recall that ${\cal S}_{\be}$ is spanned by $D_{p,q}(x^{\be,\vec i})$ with $\vec i\in\NJ$ and
$p,q\in\ol{1,\ell}$. Using induction on $\vec j$, one can prove that the
element $d$ is a linear combination of
$$D_{p,q}(x^{\be,\vec i})\;\mbox{ with }\;p,q\in\ol{1,\ell_1}\;\mbox{ or }\;p\in\ol{1,\ell_1},\,q\in\ol{\ell_1+1,\ell}\;
\mbox{ but }\;i_p=0.\eqno(4.28)$$
Thus we can assume that $d$ has the form (4.28). If $p,q\in\ol{1,\ell_1}$, then
$$ud=uD_{p,q}(x^{\be,\vec i})=D_{p,q}(x^{\al+\be,\vec i})\in {\cal S}.\eqno(4.29)$$
If $p\in\ol{1,\ell_1},\,q\in\ol{\ell_1+1,\ell}$ and $i_p=0$, then
$$D_{p,q}(x^{\be,\vec i})=(\be_{q-\ell_1}-\rho_{q-\ell_1})x^{\be,\vec i}\ptl_p+i_qx^{\be,\vec i-1_{[q]}}\ptl_p.\eqno(4.30)$$
Assume $\al+\be\ne\rho$. Choose $q'\in\ol{\ell_1+1,\ell}$ with
$a=\al_{q'-\ell_1}+\be_{q'-\ell_1}-\rho_{q'-\ell_1}\ne0$. Then
\begin{eqnarray*}& &ud-(\be_{q-\ell_1}-\rho_{q-\ell_1})a^{-1}D_{p,q'}(x^{\al+\be,\vec i})\\ &=&i_qx^{\al+\be,\vec i-1_{[q]}}\ptl_p-(\be_{q-\ell_1}-\rho_{q-\ell_1})a^{-1}i_{q'}x^{\al+\be,\vec i-1_{[q']}}\ptl_p.\hspace{5.1cm}(4.31)\end{eqnarray*}
by induction on $\vec i$, we prove (ii). The case $\al+\be=\rho$ and $\ell_2\ge1$ can be proved similarly.$\qquad\Box$
\psp

{\bf Theorem 4.4 (Main Theorem)}. {\it The Lie algebras} ${\cal S}={\cal S}(\ell_1,\ell_2,\ell_3;\rho,\G)$ {\it with} $\ell\geq 3$ {\it and} ${\cal S}'={\cal S}(\ell_1',\ell_2',\ell_3';\rho',\G')$ {\it (cf. (1.26)) are isomorphic if and only if} $(\ell_1,\ell_2,\ell_3)=(\ell_1',\ell_2',\ell_3')$ {\it and there exists an element} $g\in  G_{\ell_2,\ell_3}$ {\it (cf. (1.18)) such that} $g(\G)=\G'$, {\it and} $g(\rho)=\rho'$ {\it (cf. (1.19)) if} $\ell_1=0$. {\it In particular, there exists a one-to-one correspondence between the set of isomorphism classes of the Lie algebras of the form (1.26) and the set} $SW$ {\it in (1.23) if} $\ell_1>0$, {\it and between the set of isomorphism classes of the Lie algebras of the form (1.26) and the set} $SS$ {\it defined in (1.28)  if} $\ell_1=0$.
\psp

{\it Proof.} We shall use the notations with a prime to denote elements and vector spaces related to ${\cal S}'$; for instance, ${\cal A}'$, ${x'}^{\al',\vec{i}}$. Let $\sgm: {\cal S}\rar {\cal S}'$ be a Lie algebra isomorphism. Then $\sgm$ induces an isomorphism, also denoted by $\sgm$, of $\der{\cal S}$ to $\der{\cal S'}$. For convenience, we view ${\cal S}$ as a Lie
subalgebra of $\der{\cal S}$ by identifying $u\in {\cal S}$ with $\ad_u$ when the context is clear.
\pse

\ul{\it Case 1}. $\ell_1=0$.
\pse

By Lemma 4.1, the maximal space of locally-finite inner derivations of
${\cal S}$ is ${\cal D}_{\rho}$ which is of finite dimensional. So any maximal space of locally-finite inner derivations of ${\cal S}'$ is of finite dimensional. Again by Lemma 4.1, we must have
$$\ell'_1=0=\ell_1\mbox{ \ and \ }\dim{\cal D}_{\rho}=\dim{\cal D}'_{\rho'}.
\eqno(4.32)$$
If $\ell_2=0$, then by Lemma 4.1, ${\cal S}$ has no locally-nilpotent derivations and so does ${\cal S}'$. Hence  $\ell'_2=0$ by (4.3). Our theorem follows from Theorem 5.1 in [DZ]. Now  we assume $\ell_2,\ell'_2\ge1$.

By Lemma 4.1, ${\rm Hom}_{\BB{Z}}(\G,\BB{F})$,
${\cal D}_3\cap{\cal D}_{\rho}$ and ${\cal D}_{\rho}$
are the space of semi-simple derivations, the space of  semi-simple
inner derivations and the space of locally-finite inner derivations of ${\cal S}$, respectively. Thus
$$\sgm({\rm Hom}_{\BB{Z}}(\G,\BB{F}))={\rm Hom}_{\BB{Z}}(\G',\BB{F}),\ \ \ \sgm({\cal D}_3\cap{\cal D}_{\rho})={\cal D}'_3\cap{\cal D}'_{\rho'},\ \ \ \sgm({\cal D}_{\rho})={\cal D}'_{\rho'}.\eqno(4.33)$$
Note that the algebra ${\cal S}=\bigoplus_{\al\in\G}{\cal S}_{\al}$ (cf. (2.9)) is a $\G$-graded Lie algebra, whose homogeneous components ${\cal S}_{\al}$ are precisely the weight spaces of ${\rm Hom}_{\BB{Z}}(\G,\BB{F})$. Thus there exists an additive group isomorphism
$$g:\G\rar\G'\;\mbox{ with }\;g(0)=0\;\mbox{ such that }\;\sgm({\cal S}_{\al})={\cal S}'_{g(\al)}\ \ \for\;\,\al\in\G.\eqno(4.34)$$
In particular,
$${\cal S}_0\cong {\cal S}'_0,\ \ (\der{\cal S})_0\cong(\der{\cal S}')_0.\eqno(4.35)$$

We want to prove that for any $\al\in\G$ and $\ptl\in{\cal D}_{\al-\rho}$,
there exists a unique $\ptl'\in{\cal D}'_{\tau(\al)-\rho'}$ such that
$$\sgm(x^{\al}\ptl)={x'}^{\tau(\al)}\ptl'.\eqno(4.36)$$

Suppose that a term ${x'}^{\tau(\al),\vec i}\ptl''$ appears in $\sgm(x^{\al}\ptl)$ with $i_q\ne0$ for some
$q\in\ol{1,\ell'_1+\ell'_2}$. Then
$$[\ptl_{t'_q},\sgm(x^{\al}\ptl)]\ne0.\eqno(4.37)$$
Since $\ptl_{t'_q}$ is locally-nilpotent derivation on ${\cal S}'$, $\sgm^{-1}(\ptl_{t'_q})$ is also locally-nilpotent derivation on ${\cal S}$. By Lemma 4.1, $\sgm^{-1}(\ptl_{t'_q})\in{\cal D}^-$. So we have
$[\sgm^{-1}(\ptl_{t'_q}),x^{\al}\ptl]=0,$ which contradicts (4.37). Thus (4.36) holds, and $\sgm$ induces an isomorphism $\sgm:{\cal S}^{[0]}\cong {\cal S}'{}^{[0]}$.  By Theorem 5.1 in [DZ],
$\sgm$ induces an isomorphism $\phi=\sgm|_{\cal D}:{\cal D}\rar{\cal D}'$ such that
$$g(\rho)=\rho',\;\;\la\ptl,\al\ra=\la\phi(\ptl),g(\al)\ra\,\; \for\;\;\al\in\G.\eqno(4.38)$$
In particular, $\ell_2+\ell_3=\ell'_2+\ell'_3$. Since ${\cal D}_3$ is a subspace of semi-simple derivations on ${\cal S}$ and ${\cal D}'_3$ is the subspace of all
semi-simple derivations of ${\cal D}'$, we must have $\phi({\cal D}_3)={\cal D}'_3$,
and $\ell_3=\ell'_3$. Moreover, $g\in G_{\ell_2,\ell_3}$ (cf. (1.18) and (1.19)) by the second equation in (4.38).
\pse

\ul{\it Case 2}. $\ell_1\ge2$.
\pse

By Case 1, $\ell'_1\ge1$.
By (2.34), we can assume $\rho=\rho'=0.$ Then by (2.13), ${\cal D}={\cal D}_0\in {\cal S}$.
For $p\in\ol{1,\ell_1}$, since $\ptl_p\in {\cal S}$ is a locally-nilpotent inner derivation of ${\cal S}$,
Lemma 4.1 implies  $\sgm(\ptl_p)\in {\cal S}'\cap {\cal A}'{\cal D}'_1$. Write
$$\sgm(\ptl_p)=\sum_{\al'\in\G'_p}u'_{\al'},\qquad\mbox{where}\ \ u'_{\al'}\in {\cal A}'_{\al'}{\cal D}'_1\eqno(4.39)$$
and $\G'_p$ is a  finite subset of $\G'$. By Lemma 4.3 (ii), we have
$${x'}^{\be'}\sgm(\ptl_1)\in {\cal S}'\qquad \for\,\be'\in\G'\bs(-\G'_p),\eqno(4.40)$$
where $-\G'_p=\{-\al'\,|\,\al'\in\G'_p\}$. Set
$$B'=\{z\in {\cal A}'\,|\,z\sgm(\ptl)\in \der{\cal S}'\,\for\,\ptl\in{\cal D}_1\}.\eqno(4.41)$$
Let
$$\ol A\mbox{ be the subalgebra of ${\cal A}'$ generated by }B'.\eqno(4.42)$$

Since $\G'_0=-\cup_{p=1}^{\ell_1}\G'_p$ is a finite subset of $\G'$ and $\G'$ is a torsion-free
group, $\G'$ is generated by $\G'\bs\G'_0$. Since $A'_1$ (cf. (4.26)) is the group algebra $\BB{F}[\G']$, (4.40) shows 
$$A'_1\subset \ol A.\eqno(4.43)$$
For any $z\in B'$ and $p\in\ol{1,\ell}$, we have
$$\sgm(\ptl_p)(z)\sgm(\ptl)=[\sgm(\ptl_p),z\sgm(\ptl)]=\sgm([\ptl_p,\sgm^{-1}(z\sgm(\ptl))])\qquad\for\;\;\ptl\in{\cal D}_1.\eqno(4.44)$$
Thus $B'$ is $\sgm({\cal D})$-invariant, and so is $\ol A$.
As $\ptl_p$ is locally-nilpotent or locally-finite on ${\cal S}$
for $p\in\ol{1,\ell_1}$ or $p\in\ol{1,\ell}$, (4.44) implies that $\sgm(\ptl_p)$ is locally-nilpotent or locally-finite on $\ol A$. By these facts, we have
$$\ol A=\bigoplus_{\be\in\BB{F}^{\ell_2+\ell_3}}\ol A_{\be},\eqno(4.45)$$
with
$$\ol A_{\be}=\{z\in {\cal A}'\,|\,(\sgm(\ptl_p)-\ol\be_p)^m(z)=0
\mbox{ for $p\in\ol{\ell_1+1,\ell}$ and some }
m\in\BB{N}\},\eqno(4.46)$$
where we have written $\be=(\ol\be_{\ell_1+1},\cdots,\ol\be_{\ell})$. 
Set
$$\ol A_{\be}^{(0)}=\{z\in\ol A_{\be}\,|\,\sgm(\ptl_p)(z)=0,\sgm(\ptl_q)(z)=\ol\be_q z,\ \
\for\,\;p\in\ol{1,\ell_1},\,q\in\ol{\ell_1+1,\ell}\}.\eqno(4.47)$$
Obviously
$$\ol A_{\be}\ne\{0\}\lra \ol A_{\be}^{(0)}\ne\{0\}.\eqno(4.48)$$            
Take
$$\ol\G=\{\be\in\BB{F}^{\ell_2+\ell_3}\,|\,\ol A_{\be}^{(0)}\ne\{0\}\}.
\eqno(4.49)$$
Since $B'$ is $\sgm({\cal D})$-invariant, we have
$$B'=\bigoplus_{\be\in\ol\G}B'_{\be},\;B'_{\be}\ne\{0\}\lra {B'}_{\be}^{(0)}\ne\{0\},\; B'_{\be}=B'\cap \ol A_{\be},\;{B'}^{(0)}_{\be}=B'\cap \ol A^{(0)}_{\be}.\eqno(4.50)$$
Set
$$\ol\G'=\{\be\in\ol\G\,|\,{B'}^{(0)}_{\be}\ne\{0\}\}.\eqno(4.51)$$
Since $\ol A$ is generated by $B'$ and
$${B'}^{(0)}_{\be}{B'}^{(0)}_{\gm}\subseteq
\ol A^{(0)}_{\be+\gm},\ \ \ol A_{\al}=\sum_{\be+\gm=\al}B'_{\be}B'_{\gm}\ \ \ \for\;\,\be,\gm\in\ol\G',\,\al\in\ol\G,\eqno(4.52)$$
we have
$$\ol\G=\mbox{ the group generated by }\ol\G'.\eqno(4.53)$$

For any $\ptl\in{\cal D}_1\bs\{0\}$, there exists $z\in B'$ such that $z\sgm(\ptl)\notin\sgm({\cal D})$ because  
$${\rm Span}\{{x'}^{\al'}\sgm(\ptl)\,|\,\al'\in\G'\bs(-\G'_0)\}\;\;\mbox{is infinite-dimensional.}\eqno(4.54)$$
Write
$$z\sgm(\ptl)=\sgm(w),\ w=\sum_{\al\in G}x^{\al}u_{\al},\,u_{\al}=\sum_{\vec i\in J_{\al}}t^{\vec i}\ptl^{\al,\vec i},\eqno(4.55)$$
where $G$ is the finite subset of the elements $\al\in\G$ such that $u_{\al}\ne0$ and $J_{\al}$ is the finite subset of the elements $\vec i\in\NJ$ such that $\ptl^{\al,\vec i}\in{\cal D}\bs\{0\}$. Moreover,  $G\ne\{0\}$ or $J_0\ne\{0\}$. We claim that
$$w\in {\cal A}{\cal D}_1.\eqno(4.56)$$
If not, by Lemma 4.2, there exists $u\in{\rm Nil}({\cal S})$ such that $\ad_w$ is not nilpotent acting on $u$. Since $\sgm(u)$ is locally-nilpotent, $\sgm(u)\in {\cal S}'\cap {\cal A}'{\cal D}'_1$  by Lemma 4.1. But by Lemma 4.3 (i), $\ad_{z\sgm(\ptl_1)}$ is locally-nilpotent on ${\cal A}'{\cal D}'_1$, in particular, on $\sgm(u)$. This leads a contradiction.

Note that for any basis $\{d_1,\cdots,d_{\ell_1}\}$ of ${\cal D}_1$, there exists
an automorphism $\iota$ of ${\cal S}$ such that
$$\iota(d_p)=\ptl_p\ \ \for\;\,p\in\ol{1,\ell_1}.\eqno(4.57)$$
Assume  $G=\{0\}$ in (4.55). Then there exists $\vec i\in J_0\bs\{0\}$, say, $i_p\ne0$ with $p\in\ol{1,\ell_1+\ell_2}$.
By (4.57), we can assume $\ptl=\ptl_1$ in (4.55). Then for any $\be\in\G$, we can either take $u=x^{\be}\ptl_p$ if $\rho_{p-\ell_1}=0$, or $u=\rho_{p-\ell_1}x^{\be,1_{[2]}}\ptl_2-x^{\be}\ptl_p=
D_{2,p}(x^{\be,1_{[2]}}) \in {\cal S}$ if $\rho_{p-\ell_1}\ne0$, such that
$$\sgm(u)(z)\sgm(\ptl)=[\sgm(u),z\sgm(\ptl)]=\sgm([u,w])\ne0.\eqno(4.58)$$
Replacing $z$ by $\sgm(u)(z)$ in (4.55), we can assume $G\ne\{0\}$.
Note that for any $\td\ptl\in{\cal D}$,  we have
$$\sgm(\td\ptl)^m(z)\sgm(\ptl)=\ad_{\sgm(\td\ptl)}^m (z\sgm(\ptl))
=\sgm(\sum_{\al\in G}x^{\al}\sum_{k=0}^m\left(\!\!\!\begin{array}{c}m\\ k\end{array}\!\!\!\right)\al(\td\ptl)^k{\td\ptl}^{m-k}(u_{\al})).\eqno(4.59)$$
As $\td\ptl$ is nilpotent on $u_{\al}$, we obtain that there exists some
$z\in B'\bs\{0\}$ such that
$$z\sgm(\ptl)=\sgm(x^{\al}d)\mbox{ for some $x^{\al}d\in {\cal S}$ with }\al\in\G\bs\{0\},d\in{\cal D}_1\bs\{0\}\eqno(4.60)$$
 by induction on the leading degree of $w$. Thus
$$\sgm(\td\ptl)(z)\sgm(\ptl)=[\sgm(\td\ptl),z\sgm(\ptl)]=\sgm([\td\ptl,x^{\al}d])=\al(\td\ptl)z\sgm(\ptl).\eqno(4.61)$$
which implies $z\in {B'}_{\al}^{(0)}$. So $\al\in\ol\G'$.

We claim that for any $\be\in\G$, there exists $z'\in {B'}^{(0)}_{\be}\bs\{0\}$ such that
$$z\sgm(\ptl)=\sgm(x^{\be}d)\;\mbox{ for some }\;d\in{\cal D}_1\bs\{0\}.\eqno(4.62)$$

Take $\ptl=\ptl_1$. First suppose $d\notin\BB{F}\ptl_1$ in (4.60), say, $\chi_p(d)=1$ for some
$p\in\ol{2,\ell_1}$ (cf. (2.26)). Then for any $\be\in\G$, we have
$$\sgm(x^{\be-\al,1_{[p]}}\ptl_1)(z)\sgm(\ptl)=[\sgm(x^{\be-\al,1_{[p]}}\ptl_1),z\sgm(\ptl)]=\sgm([x^{\be-\al,1_{[p]}}\ptl_1,x^{\al}d])=-\sgm(x^{\be}\ptl_1),\eqno(4.63)$$
where the last equality follows from the fact that $d\in{\cal D}_1$; that is, (4.62) holds.
Next suppose $d=\ptl_1$ in (4.60). We choose $p$ with $\al_{p-\ell_1}\ne0$. For any $\be\in\G$ with $a=\be_{p-\ell_1}-\al_{p-\ell_1}\ne0$, as in (4.63), we have
\begin{eqnarray*}\sgm(a x^{\be-\al,1_{[2]}}\ptl_2-x^{\be-\al}\ptl_p)(z)\sgm(\ptl)&=&[\sgm(a x^{\be-\al,1_{[2]}}\ptl_2-x^{\be-\al}\ptl_p),z\sgm(\ptl)]\\ &=&-\al_{p-\ell_1}\sgm(x^{\be}\ptl_1).\hspace{5.5cm}(4.64)\end{eqnarray*}
Since (4.64) holds for all $\be$ with $\be_{p-\ell_1}\ne\al_{p-\ell_1}$, we can derive (4.62)  from (4.64).
This proves (4.62), and (4.62) proves $\G\subset\ol\G'.$ Conversely, for any $\be\in\ol\G'$ and $z\in {B'}_{\be}^{(0)}$, (4.59) implies $\be\in\G$, which together with (4.53) proves
$$\G=\ol\G'=\ol\G.\eqno(4.65)$$
By (4.62), for any $z\in {B'}_{\al}^{(0)}\bs\{0\}$ and $\al\in\ol\G'$, we can write
$$z\sgm(\ptl)=\sgm(x^{\al}\tau_z(\ptl))\;\mbox{ with }\;\tau_z(\ptl)\in{\cal D}_1\;\;\for\,\;\ptl\in{\cal D}_1.\eqno(4.66)$$
Hence we get an injective linear transformation $\tau_z$ on ${\cal D}_1$ because
${\cal A}'$ has no zero divisors. Since ${\cal D}_1$ is finite dimensional, $\tau_z$
is a linear automorphism. For any $\al\in\G$ and $z\in {B'}_{\al}^{(0)}\bs\{0\}$, we claim
that $\tau_z$ is a scalar operator. Otherwise by (4.57), we can assume
$$z\sgm(\ptl_1)=\sgm(x^{\al}\ptl_2).\eqno(4.67)$$
Take an eigenvector $d$ of $\tau_z$. By (4.67), we see that
$\ptl_1,\ptl_2,d$ are linearly independent since ${\cal A}'$ has
no zero divisors. So we can assume $d=\ptl_3$, that is,
$$z\sgm(\ptl_3)=\sgm(x^{\al}\ptl_3). \eqno(4.68)$$
Applying $\sgm(t_2\ptl_1)$ to (4.67) and (4.68), we obtain that
$\sgm(t_2\ptl_1)(z)$ is nonzero and zero, respectively. This leads a contradiction.
Thus for any $\al\in\G$ and $z\in {B'}^{(0)}_{\al}\bs\{0\}$, there exists $c_z\in\BB{F}\bs\{0\}$ such that
$$z\sgm(\ptl)=c_z\sgm(x^{\al}\ptl)\qquad\for\,\;\ptl\in{\cal D}_1.\eqno(4.69)$$
In particular,
$$\dim{B'}_{\al}^{(0)}=1\qquad \for\,\;\al\in\G,\eqno(4.70)$$
and
$$\ol A_{\al}={B'}_{\al}^{(0)}\ol A_0\qquad \for\,\;\al\in\G.\eqno(4.71)$$

For given $\al\in\G$, we denote
$$\td\ptl_1=x^{\al}\ptl_1,\,\;\td\ptl_2=x^{-\al}\ptl_2,\,\;\td\ptl_p=\ptl_p,\,\;\td\ptl_q=\ptl_q+\al_q(t_1\ptl_1-t_2\ptl_2)
\eqno(4.72)$$
for $p\in\ol{3,\ell_1}$ and $q\in\ol{\ell_1\!+\!1,\ell}$.
By (2.31) and (2.32), there exists an automorphism of ${\cal S}$ which maps $\ptl_p$ to $\td\ptl_p$ for $p\in\ol{1,\ell}$ and fixes
$A_1$. Thus we can replace $\ptl$ by $\td\ptl$ in the
above arguments. In particular, for any nonzero $z\in {B'}^{(0)}_{\al}$, we let $\td z\in \td {B'}^{(0)}_{-\al}\bs\{0\}$, where $\td {B'}^{(0)}_{-\al}$ is defined as in (4.41) and (4.50) with respect to
$\td{\cal D}_1=\{\td\ptl_p\,|\,p\in\ol{1,\ell_1}\}$. Moreover, we have
$$\td z z\sgm(\ptl_1)=\td z c_z\sgm(x^{\al}\ptl_1)=c_z\td z\sgm(\td\ptl_1)
=c_z c_{\td z}\sgm(x^{-\al}\td\ptl_1)
=c_zc_{\td z}\sgm(\ptl_1)\qquad \for\;\,\ptl\in{\cal D}_1.\eqno(4.73)$$
Thus $\td z z=c_zc_{\tilde z}\in\BB{F}\bs\{0\}$. This shows that $z$ is invertible. By the proof of Theorem 2.1 in [SXZ],
$$(\bigcup_{\al\in\G} {B'}^{(0)}_{\al})\bs\{0\}=\mbox{ the set of all invertible elements in }B'.\eqno(4.74)$$
Note that
$$\mbox{the set of all invertible elements in }{\cal A}'=(\bigcup_{\al'\in\G'}\BB{F} {x'}^{\al'})\bs\{0\}.\eqno(4.75)$$
Since all but a finite number of ${x'}^{\al'}$ are in $B'$, the set in (4.74) must be equal to the set in (4.75); equivalently, there exists  a group isomorphism $g:\G\rar\G'$
such that
$${B'}^{(0)}_{\al}=\BB{F} {x'}^{\tau(\al)}\qquad \for\,\;\al\in\G.\eqno(4.76)$$
Thus
$${B'}^{(0)}=\bigoplus_{\al\in\G}{B'}_{\al}^{(0)}\eqno(4.77)$$
is a group algebra. Furthermore, by (4.71),
$$\ol A={B'}^{(0)}\ol A_{0}\cong {B'}^{(0)}\otimes\ol A_0.\eqno(4.78)$$
By (4.76), we can assume
$${x'}^{g(\al)}\sgm(\ptl)=\sgm(x^{\al}\ptl)\;\;\for\;\,\al\in\G,\,\ptl\in{\cal D}_1.\eqno(4.79)$$

For any $\ptl\in{\cal D}_2+{\cal D}_3$,  $\sgm(\ptl)\in{\cal A}'{\cal D}'_1+{\cal D}'$. Define $\phi:{\cal D}_2+{\cal D}_3\rar{\cal D}'_2+{\cal D}'_3$ by
$$\phi(\ptl)=\mbox{ the image of $\sgm(\ptl)$ under the projection: }{\cal A}'{\cal D}'_1\oplus({\cal D}'_2+{\cal D}'_3)\rar ({\cal D}'_2+{\cal D}'_3).\eqno(4.80)$$
Then by definition,
$$\la\phi(\ptl),\al'\ra {x'}^{\al'}\sgm(\ptl_1)=[\sgm(\ptl),{x'}^{\al'}\sgm(\ptl_1)]=\sgm(\ptl)({x'}^{\al'})\sgm(\ptl_1)
=\al(\ptl){x'}^{\al'}\sgm(\ptl_1)\eqno(4.81)$$
for $\al\in\G$ and $\al'=\tau(\al)$, where the first equality follows from (4.80) and the last follows from (4.76).
Thus we have the second equation of (4.38), by which  $\phi$ must be a bijection because $g:\G\rar\G'$ is an isomorphism and $\G,\G'$ are nondegenerate. In particular,
$$
\ell_2+\ell_3=\ell'_2+\ell'_3.
\eqno(4.82)$$
Since an element $\ptl\in{\cal D}_3$ is semi-simple and the elements of ${\cal D}'_2$ are not semi-simple, we must have
$\sgm(\ptl)\in {\cal A}'{\cal D}'_1+{\cal D}'_3$. Otherwise, $\sgm(\ptl)$ is not semi-simple by the proof of
Lemma 4.1. Hence
$$
\phi({\cal D}_3)\subset{\cal D}'_3,
\mbox{ which implies }\ell_3\le\ell_3'.
\eqno(4.83)$$
\par
Set
$$
{B'}^{(1)}_0=\{z\in B'_0\,|\,\sgm({\cal D})(v)\subset\BB{F}\}\subset\ol A_0.
\eqno(4.84)$$
By the proof of Theorem 2.1 in [SXZ], (4.77) and (4.78), we have
\begin{eqnarray*}\dim({B'}^{(1)}_0/\BB{F})&\leq&\mbox{the transcendental degree of $\ol A_0$ over }\BB{F}\\ &\leq &\mbox{the transcendental degree of $A'_0$ over }\BB{F}=\ell'_1+\ell_2'\hspace{2.9cm}(4.85)\end{eqnarray*}
(cf. (4.26)). We claim that for any $q\in\ol{1,\ell_1+\ell_2}$, there exists $z_q\in {B'}^{(1)}_0$ such that
$$z_1\sgm(\ptl_2)=\sgm(t_1\ptl_2),\;\;z_q\sgm(\ptl_1)=\sgm(t_q\ptl_1)\qquad\for\;\; q\in\ol{2,\ell_1+\ell_2}.\eqno(4.86)$$
Assume that (4.86) holds. Applying $\sgm(\ptl_p)$, $p\in\ol{1,\ell}$ to (4.86), we obtain
$$\sgm(\ptl_p)(z_q)=\dlt_{p,q}\;\; \for\;\,p\in\ol{1,\ell},\,q\in\ol{1,\ell_1+\ell_2},\eqno(4.87)$$
which implies that $\{1,z_q\,|\,q\in\ol{1,\ell_1+\ell_2}\}$ is an linearly independent subset of ${B'}^{(1)}_0$.
Thus $\ell_1+\ell_2\leq(\dim{B'}^{(1)}_0/\BB{F})\leq\ell'_1+\ell'_2$.
This together with (4.82), (4.83) gives $\ell'_1\ge\ell_1\ge2$. Hence
we can exchange positions of ${\cal S}$ and ${\cal S}'$ so that we obtain
$\ell'_3\le\ell_3$ and $\ell'_1+\ell'_2\le\ell_1+\ell_2$,
which together with (4.82) and (4.83) imply our theorem in this case.

Now we want to prove (4.86). Choose $\al\in\G\bs\{0\}$ with $\al_p\ne0$ for some $p$. Let $z\in {B'}_{\al}^{(0)}$ such that
$$z\sgm(\ptl_1)=\sgm(x^{\al}\ptl_1).\eqno(4.88)$$
If $q\in\ol{2,\ell_1}$ or $\al_{q-\ell_1}=0$, then $\al_{p-\ell_1}x^{-\al,2_{[q]}}\ptl_q+2x^{-\al,1_{[q]}}\ptl_p=D_{p,q}(x^{-\al,2_{[q]}})\in {\cal S}$. Moreover, we have
\begin{eqnarray*}\hspace{2cm}& &\sgm(\al_{p-\ell_1}x^{-\al,2_{[q]}}\ptl_q+2x^{-\al,1_{[q]}}\ptl_p)(z)\sgm(\ptl_1)\\
&=&[\sgm(\al_{p-\ell_1}x^{-\al,2_{[q]}}\ptl_q+2x^{-\al,1_{[q]}}\ptl_p),z\sgm(\ptl_1)]\\ &=&\sgm([\al_{p-\ell_1}x^{-\al,2_{[q]}}\ptl_q+2x^{-\al,1_{[q]}}\ptl_p,x^{\al}\ptl_1])\\ &=&2\al_{p-\ell_1}\sgm(t_q\ptl_1).\hspace{9.4cm}(4.89)\end{eqnarray*}
Thus we have (4.86). If $q\in\ol{\ell_1+1,\ell_1+\ell_2}$ and $\al_q\ne0$, as in (4.89), we have
$$\sgm(\al_{q-\ell_1}x^{-\al,1_{[2]}+1_{[q]}}\ptl_2-x^{-\al,1_{[2]}}\ptl_2+
x^{-\al,1_{[q]}}\ptl_q)(z)\sgm(\ptl_1)=\al_{q-\ell_1}\sgm(t_q\ptl_1).\eqno(4.90)$$
The first equation of (4.86) can be obtained from the second by exchanging positions of $q$ and $1$.
\pse

\ul{\it Case 3}. $\ell_1=1$.
\pse

 By Cases 1 and 2, $\ell'_1=1$. By Lemma 2.3, we can assume $\rho\ne0\ne\rho'$. Then by (2.13), $\ptl_1\in{\cal D}_{\rho}\subset {\cal S}$. By Lemma 4.1, 
$$\{x^{\al,\vec i}\ptl_1\,|\,(\al,\vec i)\in\G\times\NJ,\,i_1=0\}\cup{\cal D}^-\eqno(4.91)$$
 spans the space of locally-nilpotent derivations of ${\cal S}$. When $\ell_2=0$, $x^{\rho}\ptl_1\notin {\cal S}$, which is
a derivation of ${\cal S}$ by Lemma 3.3. Since $\ptl_1$ is an inner derivation
of ${\cal S}$, $\sgm(\ptl_1)$ is an inner derivation of ${\cal S}'$. Thus
$$\sgm(\ptl_1)=\sum_{(\al',\vec i)\in\G'_1\times J'_1}a_{\al',\vec i}{x'}^{\al',\vec i}\ptl_1,\eqno(4.92)$$
where $a_{\al',\vec i}\in\BB{F}$ and $\G'_1\times J'_1$ is a finite subset of $\G'\times\BB{N}^{\ell_1'+\ell_2'}$ such that $i_1=0$ for all $\vec i\in J'_1$. We have
$${x'}^{\be'}\sgm(\ptl_1)\in \der{\cal S}'\qquad \for\,\;\be'\in\G'\eqno(4.93)$$
by Lemma 4.3 (ii) if $\ell'_2\ge1$. If $\ell'_2=0$, then $J'_1=\{0\}$. When $\be'+\al'=\rho'$ for some $\al'\in\G'_1$, 
${x'}^{\be'}\cdot {x'}^{\al'}\ptl'_1\in {{\cal W}'}^{[0]}_{\rho}\subset\der{\cal S}'$. By the arguments analogous to those given in Case 2 after (4.40), we can prove the theorem in this case.

This completes the proof of the theorem.$\qquad\Box$

\vspace{0.7cm}

\noindent{\Large \bf References}

\hspace{0.3cm}

\begin{description}

\item[{[Bo]}] R. E. Borcherds, Vertex algebras, Kac-Moody algebras, and the Monster,
{\it Proc. Natl. Acad. Sci. USA} {\bf 83} (1986), 3068-3071.

\item[{[DZ]}] D. \v{Z}. Djokovic and K. Zhao, Generalized Cartan type S Lie algebras in characteristic 0, {\it J. Algebra} {\bf 193} (1997), 144-179.

\item[{[FLM]}] I. B. Frenkel, J. Lepowsky and A. Meurman, {\it Vertex Operator
Algebras and the Monster}, Pure and Applied Math. Academic Press, 1988.

\item[{[K1]}] V. G. Kac, A description of filtered Lie algebras whose associated graded Lie algebras are of Cartan types, {\it Math. of USSR-Izvestijia} {\bf 8} (1974), 801-835.

\item[{[K2]}] V. G. Kac, Lie superalgebras, {\it Adv. Math.} {\bf 26} (1977), 8-96.

\item[{[K3]}] V. G. Kac, {\it Vertex algebras for beginners}, University lectures series, Vol {\bf 10}, AMS. Providence RI, 1996.

\item[{[K4]}] V. G. Kac, Classification of infinite-dimensional simple linearly compact Lie superalgebras, {\it Adv. Math.} {\bf 139} (1998), 1-55.

\item[{[K]}] N. Kawamoto, Generalizations of Witt algebras over a field of characteristic zero,
{\it Hiroshima Math. J.} {\bf 16} (1986), 417-462.

\item[{[O]}] J. Marshall Osborn, New simple infinite-dimensional Lie algebras of characteristic 0, {\it J. Algebra} {\bf 185} (1996), 820-835.

\item[{[P]}] D. P. Passman, Simple Lie algebras of Witt type, {\it J. Algebra} {\bf 206} (1998), 682-692.

\item[{[SXZ]}] Y. Su, X. Xu and H. Zhang, Derivation-simple algebras and the structures of  Lie algebras of Witt type, {\it J. Algebra}, to appear. 

\item[{[X1]}] X. Xu, New generalized simple Lie algebras of Cartan type over a field with characteristic 0, {\it J. Algebra} {\bf 224} (2000), 23-58.

\item[{[X2]}] X. Xu, {\it Introduction to Vertex Operator Superalgebras and Their Modules}, Kluwer Academic Publishers, Dordrecht/Boston/London, 1998.

\item[{[Z]}] K. Zhao, Generalized Cartan type S Lie algebras in characteristic zero (II), {\it Pacific J. Math.} {\bf 192} (2000), 431-454.

\end{description}
\end{document}